\documentstyle[11pt]{article}
\newlength{\standardunitlength}
\newtheorem{cor}{Corollary} \newtheorem{lemma}{Lemma}
\newtheorem{theorem}{Theorem} \newtheorem{prop}{Proposition}
\newenvironment{proof}{\noindent {\sc Proof:}}{$\Box$ \vspace{2 ex}}
\begin{document}

\begin{center}
Probabilistic Measures and Algorithms Arising from the Macdonald
Symmetric Functions
\end{center}

\begin{center}
By Jason Fulman
\end{center}

\begin{center}
Dartmouth College
\end{center}

\begin{center}
Jason.E.Fulman@Dartmouth.Edu
\end{center}

\begin{abstract}
	The Macdonald symmetric functions are used to define measures on
the set of all partitions of all integers. Probabilistic algorithms are
given for growing partitions according to these measures. The case of
Hall-Littlewood polynomials is related to the finite classical groups, and
the corresponding algorithms simplify. The case of Schur functions leads to
a $q$-analog of Plancharel measure, and a conditioned version of the corresponding
algorithms yields generalizations of the hook walk of combinatorics.
\end{abstract}

\section{Introduction} \label{Introduction}

	The Macdonald symmetric functions are the most general class of
symmetric functions known at present. Various specializations give the
Schur functions, the Hall-Littlewood polynomials, Jack's symmetric
functions, and zonal polynomials. A good
account of symmetric function theory is Macdonald's book \cite{Mac}.

	The present work consider probabilistic aspects of Macdonald's symmetric
functions. Our initial motivation
came from the theory of random matrices. Recently there has been much interest in
studying what a random element of a finite general linear group $GL(n,q)$ ``looks
like''
\cite{Nu3}, \cite{Di2}, \cite{FiH}, \cite{Go2}, \cite{Han}, \cite{St1},
\cite{St2}. Many properties of a matrix (e.g. its characteristic
polynomial, its order, the dimension of its fixed space) are functions only
of its conjugacy class. Thus a logical step in understanding a random
matrix is to first understand the conjugacy class of a random
matrix. Recall that the conjugacy classes of $GL(n,q)$ correspond to the
rational canonical form of a matrix (this is a generalization of Jordan
canonical form which works over non-algebraically closed fields--see
Chapter 6 of Herstein \cite{Her}). Rational canonical form admits the
following combinatorial description. To each monic irreducible polynomial
$\phi$ over $F_q$, a field of size $q$, associate a partition (perhaps the
trivial partition) $\lambda_{\phi}$ of some non-negative integer
$|\lambda_{\phi}|$. Let $m_{\phi}$ be the degree of $\phi$. Then the data
$\lambda_{\phi}$ represents a conjugacy class when:

\begin{enumerate}

\item $|\lambda_z| = 0$
\item $\sum_{\phi} |\lambda_{\phi}| m_{\phi} = n$

\end{enumerate}

	Fulman \cite{fulthesis} defines a measure on the set of all partitions $\lambda$
of all integers as follows. Fix $u$ such that $0<u<1$. Then pick the size
of a general linear group with probability of size $n$ equal to
$(1-u)u^n$. Next pick $\alpha$ uniformly in $GL(n,q)$ and take the
partition $\lambda_{\phi}(\alpha)$ corresponding to $\phi$ in the rational
canonical form of $\alpha$.

	Theorem \ref{main} of Section \ref{Classical} proves that these
group theoretic measures on partitions can be defined in terms of the
Hall-Littlewood symmetric functions. This, together with the identities in
Section \ref{Background}, led us to a general definition of measures which
works for the Macdonald symmetric functions. The measures defined from the
Macdonald symmetric functions can be grown probabilistically (Sections
\ref{MainAlg}, \ref{SimpleAlg},
\ref{Tableau}). One remarkable feature of these algorithms is that they
blend nicely with the algebraic structure of symmetric functions. For
instance, the algorithms can be divided into steps corresponding to each
variable $x_i$. The Pieri formula of algebraic geometry also makes an
appearance as a probabilistic transition rule.

	This paper is structured as follows. Section \ref{Notation}
collects notation which will be used freely in following sections. Section
\ref{Background} reviews identities satisfied by the Macdonald symmetric
functions. Section \ref{Defining} uses the Macdonald symmetric functions to
define measures on the set of all partitions of all integers. Section
\ref{MainAlg} gives probabilistic algorithms for growing partitions
according to the measures of Section \ref{Defining}. Section
\ref{SimpleAlg} shows that the algorithms of Section \ref{MainAlg} simplify
for the case of Hall-Littlewood polynomials. Section \ref{Classical}
considers further specializations of the measures coming from the
Hall-Littlewood polynomials, explaining the connection with the finite
classical groups. Section \ref{Tableau} develops the ``Young Tableau
Algorithm'', a simplification which works only in the case relevant to the
general linear groups. Section \ref{Weights} develops a formula for the
specialized Hall-Littlewood measures in terms of weights on the Young
lattice; this extends to the unitary groups as well. Section
\ref{Plancherel} specializes the measures of Section \ref{Defining} to the
Schur functions, leading to a $q$-analog of Plancherel measure. Section
\ref{Kerov} explains how the algorithm of Section \ref{MainAlg} is related
to Kerov's $q$-generalization of the hook walk of combinatorics. Section
\ref{Suggestions} gives suggestions for future research.

	 Most of the results of this paper are taken from Fulman's Ph.D. thesis
\cite{fulthesis} done under the guidance of Persi Diaconis. The purpose of
this paper is to emphasize symmetric function theory and combinatorics with
a minimum of group theory. A companion paper to this one is Fulman
\cite{fulalgorithm}, which applies the results of Sections
\ref{Tableau} and \ref{Weights} to prove group theoretic results about the
general linear and unitary groups.

\section{Notation} \label{Notation}

	We begin by reviewing some standard notation about partitions, as
on pages 2-5 of Macdonald $\cite{Mac}$. Let $\lambda$ be a partition of
some non-negative integer $|\lambda|$ into parts $\lambda_1 \geq \lambda_2
\geq \cdots$. Let $m_i(\lambda)$ be the number of parts of $\lambda$ of
size $i$, and let $\lambda'$ be the partition dual to $\lambda$ in the
sense that $\lambda_i' = m_i(\lambda) + m_{i+1}(\lambda) + \cdots$. Let
$n(\lambda)$ be the quantity $\sum_{i \geq 1} (i-1) \lambda_i = \sum_i
{\lambda_i' \choose 2}$.

	It is also useful to define the diagram associated to $\lambda$ as
the set of points $(i,j) \in Z^2$ such that $1 \leq j \leq \lambda_i$. We
use the convention that the row index $i$ increases as one goes downward
and the column index $j$ increases as one goes across. So the diagram of
the partition $(5441)$ is:

\[ \begin{array}{c c c c c}
		. & . & . & . & .  \\
		. & . & . & . &    \\
		. & . & . & . &    \\
		. & & & &  
	  \end{array}  \]

	It is sometimes useful to think of these dots as boxes. Given a
square $s$ in the diagram of a partition $\lambda$, let $l_{\lambda}'(s),
l_{\lambda}(s), a_{\lambda}(s), a_{\lambda}'(s)$ be the number of squares
in the diagram of $\lambda$ to the north, south, east, and west of $s$
respectively.  The subscript $\lambda$ will sometimes be omitted if the
partition $\lambda$ is clear from context. So the diagram

\[ \begin{array}{c c c c c}
		. & . & . & . & .  \\
		. & s & . & . &    \\
		. & . & . & . &    \\
		. & & & &  
	  \end{array}  \]

	has $l'(s)=l(s)=a'(s)=1$ and $a(s)=2$.

	A skew-diagram is the set theoretic difference $\lambda - \mu$ of
two diagrams $\lambda$ and $\mu$, where the diagram of $\lambda$ contains
the diagram of $\mu$. A horizontal strip is a skew-diagram with at most one
square in each column. For instance the following diagram is a horizontal
strip:

\[ \begin{array}{c c c c c}
		 &  &  &  & .  \\
		 &  &  &  &    \\
		&  & . & . &    \\
		. & & & &  
	  \end{array}  \]

	Letting $f(u)$ be a polynomial in the variable $u$, the notation
$[u^n] f(u)$ means the coefficient of $u^n$ in $f(u)$. 

	The following notation is less widely known, and is taken from
Chapter 6 of Macdonald \cite{Mac}.

\begin{enumerate}

\item Given a partition $\lambda$ and a square $s$, set $b_{\lambda}(s)=1$
if $s \not \in \lambda$. Otherwise set:

\[ b_{\lambda}(s) = \frac{1-q^{a_{\lambda}(s)} t^{l_{\lambda}(s)+1}}
{1-q^{a_{\lambda}(s)+1}t^{l_{\lambda}(s)}} \]

	Let $b_{\lambda}(q,t) = \prod_{s \in \lambda} b_{\lambda}(s)$.

\item Define

\[ \phi_{\lambda / \mu}(q,t) = \prod_{s \in C_{\lambda / \mu}}
\frac{b_{\lambda}(s)}{b_{\mu}(s)} \]

	where $C_{\lambda / \mu}$ is the union of the columns
intersecting $\lambda - \mu$.

\item The skew Macdonald polynomials (in one variable) are defined as:

\[ P_{\lambda / \mu}(x;q,t) = \frac{b_{\mu}(q,t)}{b_{\lambda}(q,t)}
\phi_{\lambda / \mu}(q,t)x^{|\lambda|-|\mu|} \]

	if $\lambda - \mu$ is a horizontal strip, and $0$ otherwise.

\item Let $(x,q)_{\infty}$ denote $\prod_{i=1}^{\infty}
(1-xq^{i-1})$. Then define $\prod (x,y;q,t)$ by:

\[ \prod (x,y;q,t) = \prod_{i,j=1}^{\infty} \frac{(tx_iy_j,q)_{\infty}}
{(x_iy_j,q)_{\infty}} \]

	Also define $g_n(y;q,t)$ as the coefficient of $x^n$ in
$\prod_j \frac{(txy_j,q)_{\infty}} {(xy_j,q)_{\infty}}$.

\end{enumerate}

\section{Properties of the Macdonald Symmetric Functions}
\label{Background}

	The Macdonald symmetric functions $P_{\lambda}(x_i;q,t)$ are a
two-parameter family of symmetric functions. A precise definition is in Chapter 6
of Macdonald \cite{Mac}. The Macdonald symmetric functions have five properties
which we shall need. It is convenient to name them (the Pieri Formula is already
named).

\begin{enumerate}

\item Measure Identity $\cite{Mac}$, page 324:

\[ \sum_{\lambda} P_{\lambda}(x;q,t) P_{\lambda}(y;q,t)
b_{\lambda}(q,t) = \prod (x,y;q,t) \]

\item Factorization Theorem $\cite{Mac}$, page 310:

\[ \prod (x,y;q,t) = \prod_{n \geq 1} e^{\frac{1}{n}
\frac{1-t^n}{1-q^n} p_n(x) p_n(y)} \]

\item Principal Specialization Formula $\cite{Mac}$, page 337:

\[ P_{\lambda}(1,t,\cdots,t^{N-1};q,t) = t^{n(\lambda)} \prod_{s
\in \lambda} \frac{1-q^{a'(s)}t^{N-l'(s)}}{1-q^{a(s)}t^{l(s)+1}} \]

\item Skew Expansion $\cite{Mac}$, pages 343-7:

\[ P_{\lambda}(x_1,\cdots,x_N;q,t) = \sum_{\mu}
P_{\mu}(x_1,\cdots,x_{N-1};q,t) P_{\lambda / \mu}(x_N;q,t) \]

\item Pieri Formula $\cite{Mac}$, page 340:

\[ P_{\mu}(y;q,t) g_r(y;q,t) = \sum_{\lambda: |\lambda - \mu|=r \atop
\lambda - \mu \ horiz. \ strip} \phi_{\lambda / \mu}(q,t)
P_{\lambda}(y;q,t) \]

	It is worth remarking that the Pieri Formula has its history in algebraic
geometry, as a rule for multiplying classes of Schubert varieties in
the cohomology ring of Grassmanians.

\end{enumerate}

\section {Defining Measures $P_{x,y,q,t}$ from the Macdonald
Symmetric Functions} \label{Defining}

	In this section the Macdonald symmetric functions are used to
define families of probability measures on the set of all partitions of all
numbers. It is assumed throughout this paper that $x,y,q,t$ satisfy the
following conditions:

\begin{enumerate}

\item $0 \leq t,q<1$
\item $x_i,y_i \geq 0$ 
\item $\sum_{i,j} \frac{x_iy_j}{1-x_iy_j} < \infty$

\end{enumerate}

	The following formula defines a probability measure
$P_{x,y,q,t}$ on the set of all partitions of all numbers:

\[ P_{x,y,q,t}(\lambda) = \frac{P_{\lambda}(x;q,t) P_{\lambda}(y;q,t)
b_{\lambda}(q,t)}{\prod (x,y;q,t)} \]

\begin{lemma} \label{measure} $P_{x,y,q,t}$ is a measure.
\end{lemma}

\begin{proof}
	By the Measure Identity and the fact that there are countably
many partitions, it suffices to check that $0 \leq P_{x,y,q,t}
(\lambda) < \infty$ for all $\lambda$. For this it is sufficient to
show (again by the Measure Identity) that $P_{\lambda}(x;q,t),
b_{\lambda}(q,t) \geq 0$ for all $\lambda$ and that $0 \leq \prod
(x,y;q,t) < \infty$.

	Condition 1 implies that $b_{\lambda}(q,t) \geq 0$ for all
$\lambda$. We claim that $x_i \geq 0$ implies that $P_{\lambda}(x;q,t)
\geq 0$. To see this, note that when $P_{\lambda}(x;q,t)$ is expanded
in monomials in the $x$ variables, all coefficients are
non-negative. For any particular monomial, this follows by repeated
use of the Skew Expansion.

	By the Factorization Theorem, showing that $0 \leq \prod
(x,y;q,t) < \infty$ is equivalent to showing that:

\[ 0 \leq \sum_{n \geq 1} \frac{1}{n} \frac{1-t^n}{1-q^n} p_n(x)
p_n(y) < \infty \]

	Conditions 1 and 2 imply that this expression is
non-negative. To see that it is finite, use Condition 3 as follows:

\begin{eqnarray*}
\sum_{n \geq 1} \frac{1}{n} \frac{1-t^n}{1-q^n} p_n(x) p_n(y) & \leq & \frac{1}{1-q} \sum_{n \geq 1} p_n(x) p_n(y)\\
& = & \frac{1}{1-q} \sum_{i,j \geq 1} \frac{x_i y_j}{1-x_iy_j}\\
& < & \infty
\end{eqnarray*}

\end{proof}

	Define truncated measures $P^N_{x,y,q,t}(\lambda)$ to be $0$
if $\lambda$ has more than $N$ parts, and otherwise:

\[ P^N_{x,y,q,t}(\lambda) = \frac{ P_{\lambda}
(x_1,\cdots,x_N,0,\cdots;q,t) P_{\lambda}(y;q,t) b_{\lambda}(q,t)}{\prod
(x_1,\cdots,x_N,0,\cdots,y;q,t)} \]

	Let $P^0(x,y,q,t)$ be 1 on the empty partition and 0
elsewhere. Arguing as in Lemma $\ref{measure}$ shows that the
$P^N_{x,y,q,t}$ are probability measures. It is also clear that
$lim_{N \rightarrow \infty} P^N_{x,y,q,t} = P_{x,y,q,t}$. There are
other possible definitions of $P^N_{x,y,q,t}$ which converge to
$P_{x,y,q,t}$ in the $N \rightarrow \infty$ limit (for instance one
can truncate both the $x$ and $y$ variables). These deserve further
investigation.

\section{A Probabilistic Algorithm for Picking From $P_{x,y,q,t}$}
\label{MainAlg}

	This section gives a stochastic method for picking from
$P_{x,y,q,t}$ under conditions 1-3 of Section $\ref{Defining}$.

\begin{center}
Algorithm for Picking from $P_{x,y,q,t}$
\end{center}

\begin{description}

\item [Step 0] Start with $\lambda$ the empty partition and $N$ (which
we call the interval number) equal to $1$.

\item [Step 1] Pick an integer $n_N$ so that $n_N=k$ with probability
$\prod_j \frac{(x_Ny_j,q)_{\infty}}{(tx_Ny_j,q)_{\infty}} g_k(y;q,t)
x_N^k$. (These probabilities sum to 1 by the definition of $g_k$).

\item [Step 2] Let $\Lambda$ be a partition containing $\lambda$ such
that the difference $\Lambda - \lambda$ is a horizontal strip of size
$n_N$. There are at most a finite number of such $\Lambda$. Change
$\lambda$ to $\Lambda$ with probability:

\[ \frac{\phi_{\Lambda / \lambda}(q,t)}{g_{n_N}(y;q,t)}
\frac{P_{\Lambda}(y;q,t)} {P_{\lambda} (y;q,t)} \]

	(These probabilities sum to 1 by the Pieri Formula). Then set
$N=N+1$ and go to Step 1.

\end{description}

	Lemma $\ref{stop}$ will show that this algorithm terminates with probability 1.

	As an example of the algorithm, suppose we are at Step 1 with
$N=3$ and the partition $\lambda$:
	
\[ \begin{array}{c c c}
		. & . &    \\
		. &  &      \\
		 &  &    \\  
	  \end{array}  \]

	We then pick $n_3$ according to the rule in Step 1. Suppose
that $n_3=2$. We thus add a horizontal strip of size 2 to $\lambda$,
giving $\Lambda$ equal to one the following four partitions with
probability given by the rule in Step 2:

\[ \begin{array}{c c c}
		. & . &    \\
		. & . &      \\
		. &  &    \\  
	  \end{array}  \]

\[ \begin{array}{c c c}
		. & . & .   \\
		. &  &      \\
		. &  &    \\  
	  \end{array}  \]

\[ \begin{array}{c c c}
		. & . & .   \\
		. & . &      \\
		 &  &    \\  
	  \end{array}  \]

\[ \begin{array}{c c c c}
		. & . & . & .   \\
		. &  &  &    \\
		 &  &  &  \\  
	  \end{array}  \]

	We then set $N=4$ and return to Step 1. 

\begin{lemma} \label{stop} The algorithm terminates with probability
1.
\end{lemma}

\begin{proof}
	Recall the Borel-Cantelli lemmas of probability theory, which
say that if $A_N$ are events with probability $P(A_N)$ and $\sum_N
P(A_N)< \infty$, then with probability 1 only finitely many $A_N$
occur. Let $A_N$ be the event that at least one box is added to the
partition during interval $N$. To prove the lemma it is sufficient to
show that only finitely many $A_N$ occur.

	The Factorization Theorem implies that $g_0=1$. Again using the
Factorization Theorem and the fact that $1-e^{-x} \leq x$ for $x \geq
0$ shows that:

\begin{eqnarray*}
\sum_{N \geq 1} P(A_N) & = & \sum_{N \geq 1} [1-(\prod_j \frac{(x_Ny_j;q)_{\infty}}{(tx_Ny_j;q)_{\infty}}) g_0]\\
& = & \sum_{N \geq 1} [1 - e^{-\sum_{n \geq 1} \frac{1}{n} \frac{1-t^n}{1-q^n} (x_N)^n p_n(y)}]\\
& \leq & \sum_{N \geq 1} \sum_{n \geq 1} [\frac{1}{n} \frac{1-t^n}{1-q^n} (x_N)^n p_n(y)]\\
& = & \sum_{n \geq 1} \frac{1}{n} \frac{1-t^n}{1-q^n} p_n(x)p_n(y)\\
& \leq & \frac{1}{1-q} \sum_{n \geq 1} p_n(x) p_n(y)\\
& = & \frac{1}{1-q} \sum_{i,j \geq 1} \frac{x_iy_j}{1-x_iy_j}\\
& < & \infty 
\end{eqnarray*}

\end{proof}

	Theorem $\ref{mainalg}$ is the main result of this
section. Since $q$ and $t$ are fixed, the notation in the proof of
Theorem $\ref{mainalg}$ will be abbreviated somewhat by omitting the
explicit dependence on these variables.

\begin{theorem} $\label{mainalg}$ The chance that the algorithm yields
the partition $\lambda$ at the end of interval $N$ is
$P^N_{x,y,q,t}(\lambda)$. Consequently, the algorithm for picking from
$P_{x,y,q,t}$ works. \end{theorem}

\begin{proof}
	Since the algorithm proceeds by adding horizontal strips, it
is clear that the partition produced at the end of interval $N$ has at
most $N$ parts.

	The base case $N=0$ is clear since the algorithm starts with
the empty partition and $P^0_{x,y,q,t}$ is 1 on the empty partition
and 0 elsewhere.

	For the induction step, the Skew Expansion gives:

\begin{eqnarray*}
P^N_{x,y,q,t}(\Lambda) & = & [\prod_{i=1}^N \prod_j
\frac{(x_iy_j,q)_{\infty}}{(tx_iy_j,q)_{\infty}}]
P_{\Lambda}(x_1,\cdots,x_N) P_{\Lambda}(y) b_{\Lambda}\\
& = & [\prod_{i=1}^N \prod_j
\frac{(x_iy_j,q)_{\infty}}{(tx_iy_j,q)_{\infty}}] P_{\Lambda}(y)
b_{\Lambda} \sum_{\lambda \subset \Lambda}
P_{\lambda}(x_1,\cdots,x_{N-1}) P_{\Lambda / \lambda}(x_N)\\
& = & [\prod_{i=1}^N \prod_j
\frac{(x_iy_j,q)_{\infty}}{(tx_iy_j,q)_{\infty}}] P_{\Lambda}(y)
b_{\Lambda} \sum_{\lambda \subset \Lambda \atop \Lambda - \lambda \ 
horiz. \ strip} P_{\lambda}(x_1,\cdots,x_{N-1}) x_N^{|\Lambda|-|\lambda|}
\frac{b_{\lambda}}{b_{\Lambda}} \phi_{\Lambda / \lambda}\\
& = & \sum_{\lambda \subset \Lambda \atop \Lambda - \lambda \ horiz. \
strip} [(\prod_{i=1}^{N-1} \prod_j
\frac{(x_iy_j,q)_{\infty}}{(tx_iy_j,q)_{\infty}})
P_{\lambda}(x_1,\cdots,x_{N-1}) P_{\lambda}(y) b_{\lambda}]\\
&   & [ \prod_j \frac{(x_Ny_j,q)_{\infty}}{(tx_Ny_j,q)_{\infty}}
g_{|\Lambda|-|\lambda|}(y) x_N^{|\Lambda|-|\lambda|}]
[\frac{\phi_{\Lambda / \lambda}}{g_{|\Lambda|-|\lambda|}(y)}
\frac{P_{\Lambda}(y)} {P_{\lambda}(y)}]\\
& = & \sum_{\lambda \subset \Lambda \atop \Lambda - \lambda \ horiz. \
strip} [P^{N-1}_{x,y,q,t}(\lambda)] [ \prod_j
\frac{(x_Ny_j,q)_{\infty}}{(tx_Ny_j,q)_{\infty}}
g_{|\Lambda|-|\lambda|}(y) x_N^{|\Lambda|-|\lambda|}]\\
& & \ \ \ \ \ [\frac{\phi_{\Lambda / \lambda}}{g_{|\Lambda|-|\lambda|}(y)}
\frac{P_{\Lambda}(y)} {P_{\lambda}(y)}]
\end{eqnarray*}

	Probabilistically, this equality says that the chance that the
algorithm gives $\Lambda$ at the end of interval $N$ is equal to the
sum over all $\lambda$ such that $\Lambda / \lambda$ is a horizontal
strip of the chance that the algorithm gives $\lambda$ at the end of
interval $N-1$ and that $\lambda$ then grows to $\Lambda$ in interval
$N$. This proves the theorem.
\end{proof}

	As a corollary of the above algorithm, one obtains a probability
generating function with the size of the partition $\lambda$.

\begin{cor} \label{size} The distribution of the size of a partition
$\lambda$ chosen from $P_{x,y,q,t}$ has as its probability generating
function in the variable $z$:

\[ \frac{\prod (xz,y;q,t)}{\prod (x,y;q,t)} \]

\end{cor}

\begin{proof}
	By the way the algorithm works, the growth of $\lambda$ during
different intervals is independent. So it suffices to show that the chance
$\lambda$ grows by $k$ in interval $N$ is:

\[ \prod_j \frac{(x_Ny_j,q)_{\infty}}{(tx_Ny_j,q)_{\infty}} [z^k] \prod_j
\frac{(tx_Nzy_j,q)_{\infty}}{(x_Nzy_j,q)_{\infty}} \]

	This is clear from Step 1 of the algorithm and the definition
of $g_k$.
\end{proof}

	This section closes by noting that in the case $y_i=t^{i-1}$, there
is a nice expression for $g_n$. For this and future use, recall the
following lemma of Stong $\cite{St1}$.

\begin{lemma} \label{Stong} For $|q|>1$ and $0<u<1$, 
\begin{enumerate}

\item $\prod_{r=1}^{\infty} (\frac{1}{1-\frac{u}{q^r}}) = \sum_{n=0}^{\infty}
\frac{u^nq^{{n \choose 2}}}{(q^n-1) \cdots (q-1)}$

\item $\prod_{r=1}^{\infty} (1-\frac{u}{q^r}) = \sum_{n=0}^{\infty}
\frac{(-u)^n}{(q^n-1) \cdots (q-1)}$

\end{enumerate}
\end{lemma}

\begin{cor} \label{Addn} If $0<t,q<1$, then $g_n(t^{i-1};q,t) =
\frac{1}{(1-q^n) \cdots (1-q)}$

\end{cor}

\begin{proof}
	By Lemma $\ref{Stong}$,

\begin{eqnarray*}
g_n(t^{i-1};q,t) & = & [u^n] \prod_{i=1}^{\infty} (\frac{1}{1-uq^{i-1}})\\
& = & \frac{1}{q^n} [u^n]  \prod_{i=1}^{\infty} (\frac{1}{1-uq^i})\\
& = & \frac{1}{q^n} \frac{1}{q^{{n \choose 2}} (\frac{1}{q^n}-1) \cdots (\frac{1}{q}-1)}\\
& = & \frac{1}{(1-q^n) \cdots (1-q)}
\end{eqnarray*}
\end{proof}

\section{Hall-Littlewood Polynomials: Simplified Algorithms} \label{SimpleAlg}

	In this section the measure $P_{x,y,q,t}$ is studied under the
specialization $y^i=t^{i-1},q=0$. As one motivation for these choices,
note that setting $q=0$ in the Macdonald polynomials gives the
Hall-Littlewood polynomials. The further specialization $x_i=ut^i$ will be
considered in Sections \ref{Classical} - \ref{Weights}. This further
specialization is the case relevant to the finite classical groups. Nevertheless,
this section will show that the probabilistic algorithm of Section
$\ref{MainAlg}$ simplifies without having to assume that $x_i=ut^i$.

	Supposing that $0<t,x_i<1,\sum_i x_i<1$, we give a simplified
algorithm which allows one to grow the partition $\lambda$ by adding 1
box at a time. Using the Borel-Cantelli lemmas it is straightforward
to check that this algorithm always halts.

\begin{center} Simplified Algorithm for Picking from
$P_{x,t^{i-1},0,t}$
\end{center}

\begin{description}

\item [Step 0] Start with $\lambda$ the empty partition and
$N=1$. Also start with a collection of coins indexed by the natural
numbers such that coin $i$ has probability $x_i$ of heads and
probability $1-x_i$ of tails.

\item [Step 1] Flip coin $N$.

\item [Step 2a] If coin $N$ comes up tails, leave $\lambda$ unchanged,
set $N=N+1$ and go to Step 1.

\item [Step 2b] If coin $N$ comes up heads, let $j$ be the number of
the last column of $\lambda$ whose size was increased during a toss of
coin $N$ (on the first toss of coin $N$ which comes up heads, set
$j=0$). Pick an integer $S>j$ according to the rule that $S=j+1$ with
probability $t^{\lambda_{j+1}'}$ and $S=s>j+1$ with probability
$t^{\lambda_s'} - t^{\lambda_{s-1}'}$ otherwise. Then increase the
size of column $S$ of $\lambda$ by 1 and go to Step 1.

\end{description}

	For example, suppose we are at Step 1 with $\lambda$ equal to
the following partition:
	
\[ \begin{array}{c c c c}
		. & . & . & .  \\
		. & . &  &      \\
		. &  &  &    \\
		 & & &  
	  \end{array}  \]

	Suppose also that $N=4$ and that coin 4 had already come up
heads once, at which time we added to column 1, giving $\lambda$. Now
we flip coin 4 again and get heads, going to Step 2b. We have that
$j=1$. Thus we add a dot to column $1$ with probability $0$, to column
$2$ with probability $t^2$, to column $3$ with probability $t-t^2$, to
column $4$ with probability $0$, and to column $5$ with probability
$1-t$. We then return to Step 1.

	Note that the dots added during the tosses of a given coin
form a horizontal strip.

	Theorem $\ref{simplif}$ shows that the simplified algorithm works.

\begin{theorem} \label{simplif} The simplified algorithm for picking
from $P_{x,t^{i-1},0,t}$ refines the general algorithm.  \end{theorem}

\begin{proof}
	Let interval $N$ denote the time between the first and last
tosses of coin $N$. To prove the theorem, it will be shown that the
two algorithms add horizontal strips in the same way during interval
$N$.

	For this observe that the size of the strips added in interval
$N$ is the same for the two algorithms. Since $q=0$ the integer $n_N$
in Step 1 of the general algorithm is equal to $k$ with probability
$(1-x_N)x_N^k$. This is equal to the chance of $k$ heads of coin $N$
in the simplified algorithm.

	Given that a strip of size $k$ is added during interval $N$,
the general algorithm then increases $\lambda$ to $\Lambda$ with
probability:

\[ \frac{\phi_{\Lambda / \lambda}(0,t)}{g_k(t^{i-1};0,t)}
\frac{P_{\Lambda}(1,t,t^2,\cdots;0,t)} {P_{\lambda} (1,t,t^2,\cdots;0,t)}
\]

	This probability can be simplified. Lemma $\ref{Addn}$ shows
that $g_k(t^{i-1};0,t)=1$. The definition of $\phi_{\Lambda /
\lambda}(0,t)$ and the Principal Specialization Formula show that the
probability can be rewritten as:

\[ (\prod_{s \in C_{\Lambda / \lambda}} \frac{b_{\Lambda}(s)}
{b_{\lambda}(s)}) \frac{t^{n(\Lambda)} \prod_{s \in \Lambda}
\frac{1}{1-0^{a_{\Lambda}(s)}t^{l_{\Lambda}(s)+1}}} {t^{n(\lambda)}
\prod_{s \in \lambda} \frac{1}{1-0^{a_{\lambda}(s)}t^{l_{\lambda}(s)+1} }}
\]

	where $0^0=1$. Let $A$ be the set of column numbers $a>1$ such
that $\Lambda - \lambda$ intersects column $a$ but not column
$a-1$. Let $A'$ be the set of column numbers $a$ such that either
$a=1$ or $a>1$ and $\Lambda - \lambda$ intersects both columns $a$ and
$a-1$. Most of the terms in the above expression cancel, giving:

\[ \frac{t^{n(\Lambda)}}{t^{n(\lambda)}} \prod_{a \in A}
(1-t^{\lambda_{a-1}' - \lambda_a'}) = \prod_{a \in A'} t^{\lambda_a'}
\prod_{a \in A} (t^{\lambda_a'} -t^{\lambda_{a-1}'}) \]

	It is easily seen that the simplified algorithm can go from
$\lambda$ to $\Lambda$ in exactly 1 way, and this also happens with
probability equal to:

\[ \prod_{a \in A'} t^{\lambda_a'} \prod_{a \in A} (t^{\lambda_a'}
-t^{\lambda_{a-1}'}) \]

\end{proof}

\section{Hall-Littlewood Polynomials: Relation with the Finite Classical
Groups} \label{Classical}

	This section explains the relation of measures defined from the
Hall-Littlewood polynomials with the finite classical groups. The case of
the general linear groups will be worked out in detail. Analogous results
will then be described for the other classical groups.

	Recall that the conjugacy classes of $GL(n,q)$ correspond to the
possible rational canonical forms of a matrix. Rational canonical form is a
generalization of Jordan canonical form which works over non-algebraically
closed fields. See Chapter 6 of Herstein \cite{Her} for a clear discussion
of canonical forms of matrices. Rational canonical form admits the
following combinatorial description. To each monic irreducible polynomial
$\phi$ over $F_q$, a field of size $q$, associate a partition (perhaps the
trivial partition) $\lambda_{\phi}$ of some non-negative integer
$|\lambda_{\phi}|$. Let $m_{\phi}$ be the degree of $\phi$. Then the data
$\lambda_{\phi}$ represents a conjugacy class when:

\begin{enumerate}

\item $|\lambda_z| = 0$
\item $\sum_{\phi} |\lambda_{\phi}| m_{\phi} = n$

\end{enumerate}

	Given an element $\alpha$ in $GL(n,q)$, let $\lambda_{\phi}
(\alpha)$ be the partition associated to $\phi$ in the rational canonical
form of $\alpha$. For example, the identity matrix has $\lambda_{z-1}$
equal to $(1^n)$ and an elementary matrix with $a \neq 0$ in the $(1,2)$
position, ones on the diagonal and zeros elsewhere has $\lambda_{z-1}$
equal to $(2,1^{n-2})$.

	Following Fulman \cite{fulthesis}, one can now define a random partition
$\lambda_{\phi}$ as follows. Fix $u$ such that $0<u<1$. Then pick the size of a
general linear group with probability of size $n$ equal to $(1-u)u^n$. Next pick
$\alpha$ uniformly in $GL(n,q)$ and take the partition $\lambda_{\phi}(\alpha)$
corresponding to $\phi$ in the rational canonical form of $\alpha$.

	Theorem \ref{main} is the main result of this section. It shows
that the random partitions $\lambda_{\phi}$ are independent for different
$\phi$ and relates their distributions to measures defined using the
Hall-Littlewood polynomials.

\begin{theorem} \label{main} The random partitions $\lambda_{\phi}$
(defined on the union of all the groups $GL$) are independent with
distribution $P_{\frac{u^{m_{\phi}}}{q^{im_{\phi}}},
\frac{1}{q^{(i-1)m_{\phi}}},0,\frac{1}{q^{m_{\phi}}}}$.
\end{theorem}

\begin{proof}
	Kung \cite{Kun} proved that the conjugacy class of $GL(n,q)$
corresponding to the data $\lambda_{\phi}$ has size:

\[ \frac{|GL(n,q)|}{\prod_{\phi} \prod_i \prod_{k=1}^{m_i(\lambda_{\phi})}
(q^{m_{\phi}d_i(\lambda_{\phi})} - q^{m_{\phi}(d_i(\lambda_{\phi})-k)})} \]

	where for any partition $\lambda$,

\[ d_i(\lambda) = m_1(\lambda) 1 + m_2(\lambda) 2 + \cdots +
m_{i-1}(\lambda)(i-1) + (m_i(\lambda) + m_{i+1}(\lambda) + \cdots +
m_j(\lambda)) i. \]

	Next, observe that:

\begin{eqnarray*}
\prod_i \prod_{k=1}^{m_i(\lambda_{\phi})}
(q^{m_{\phi}d_i(\lambda_{\phi}} - q^{m_{\phi}(d_i(\lambda_{\phi})-k)}) & = & q^{m_{\phi}\sum_i [d_i
m_i(\lambda_{\phi}) - {m_i(\lambda_{\phi}) \choose 2}]} \prod_i
\prod_{k=1}^{m_i(\lambda_{\phi})} (1-\frac{1}{q^{km_{\phi}}})\\
& = & q^{m_{\phi} \sum_i [im_i(\lambda_{\phi})+2 \sum_{h<i} hm_h(\lambda_{\phi})]
m_i(\lambda_{\phi})} \prod_{k=1}^{m_i(\lambda_{\phi})}
(1-\frac{1}{q^{km_{\phi}}})\\
& = &  q^{m_{\phi}\sum_i (\lambda_{\phi,i}')^2} \prod_{k=1}^{m_i(\lambda_{\phi})} (1-\frac{1}{q^{km_{\phi}}})\\
& = & q^{m_{\phi}(|\lambda_{\phi}|+2n(\lambda_{\phi}))}
\prod_{k=1}^{m_i(\lambda_{\phi})} (1-\frac{1}{q^{km_{\phi}}})\\
& = & q^{m_{\phi}(|\lambda_{\phi}|+2n(\lambda_{\phi}))}
\prod_{s \in \lambda_{\phi}: a(s)=0} (1-\frac{1}{q^{(l(s)+1)m_{\phi}}})\\
& = & \frac{q^{m_{\phi} n(\lambda_{\phi})}}{P_{\lambda}(\frac{1}{q^{m_{\phi}}}
,\frac{1}{q^{2m_{\phi}}},\cdots;0,\frac{1}{q^{m_{\phi}}})}
\end{eqnarray*}

	The second equality follows from the identity $d_i(\lambda) =
[\sum_{h<i} hm_h(\lambda)] + im_i(\lambda) + [\sum_{i<k}
im_k(\lambda)]$. The third equality follows from the identity
$\lambda_i'=m_i(\lambda)+m_{i+1}(\lambda)+\cdots$. The fourth equality
follows from the identity $n(\lambda) = \sum_i {\lambda_i' \choose 2}$. The
final equality follows from the Principal Specialization Formula of Section
\ref{Background}.

	Now define a ``cycle index'' for $GL$ as in Stong \cite{St1},

\[ Z_{GL(n,q)} = \frac{1}{|GL(n,q)|} \sum_{\alpha \in GL(n,q)} \prod_{\phi
\neq z} x_{\phi,\lambda_{\phi}(\alpha)} \]

	The observation that Kung's conjugacy class size formula factors in
$\phi$ leads to the equation:

\[ 1+\sum_{n=1}^{\infty} Z_{GL(n,q)}u^n = \prod_{\phi \neq z}
\sum_{\lambda} x_{\phi,\lambda}
\frac{P_{\lambda}(\frac{u}{q^{m_{\lambda}}},\frac{u}{q^{2m_{\lambda}}} ,
\cdots;0,\frac{1}{q^{m_{\lambda}}})}{q^{m_{\phi}n(\lambda)}} \]

	The definition of the measure $P_{\frac{u^{m_{\phi}}}{q^{im_{\phi}}},
\frac{1}{q^{(i-1)m_{\phi}}},0,\frac{1}{q^{m_{\phi}}}}$ gives that:

\[ P_{\frac{u^{m_{\phi}}}{q^{im_{\phi}}},
\frac{1}{q^{(i-1)m_{\phi}}},0,\frac{1}{q^{m_{\phi}}}} (\lambda) =
\prod_{r=1}^{\infty} (1-\frac{u^{m_{\phi}}}{q^{rm_{\phi}}})
\frac{P_{\lambda}(\frac{u}{q^
{m_{\lambda}}},\frac{u}{q^{2m_{\lambda}}} ,
\cdots;0,\frac{1}{q^{m_{\lambda}}})}{q^{m_{\phi}n(\lambda)}} \]

	Therefore,

\[ 1+\sum_{n=1}^{\infty} Z_{GL(n,q)}u^n = \prod_{\phi \neq z}
\sum_{\lambda} x_{\phi,\lambda} \frac{P_{\frac{u^{m_{\phi}}}
{q^{im_{\phi}}}, \frac{1}{q^{(i-1)m_{\phi}}},0,
\frac{1}{q^{m_{\phi}}}}(\lambda)}{\prod_{r=1}^{\infty}
(1-\frac{u^{m_{\phi}}}{q^{rm_{\phi}}})} \]

	Setting all $x_{\phi,\lambda}$ to 1 in this equation gives:

\[ \frac{1}{1-u} = \prod_{\phi \neq z} \prod_{r=1}^{\infty}
\frac{1}{(1-\frac{u^{m_{\phi}}}{q^{r \ m_{\phi}}})} \]

	Combining these last two equations proves that:

\[ (1-u) [1+\sum_{n=1}^{\infty} Z_{GL(n,q)} u^n] = \prod_{\phi \neq z}
\sum_{\lambda} x_{\phi,\lambda} P_{\frac{u^{m_{\phi}}}{q^{im_{\phi}}},
\frac{1}{q^{(i-1)m_{\phi}}},0,\frac{1}{q^{m_{\phi}}}}(\lambda) \]

	The statement of the theorem is exactly a probabilistic
interpretation of this last equation.
\end{proof}	

	Theorem \ref{main} leads to a corollary which is useful for
studying the $n \rightarrow \infty$ asymptotics of random matrix theory.

 \begin{lemma} \label{bign} If $f(1)<\infty$ and $f$ has a Taylor
series around 0, then:

\[ lim_{n \rightarrow \infty} [u^n] \frac{f(u)}{1-u} = f(1) \]

\end{lemma}

\begin{proof}
	Write the Taylor expansion $f(u) = \sum_{n=0}^{\infty} a_n
u^n$. Then observe that $[u^n] \frac{f(u)}{1-u} = \sum_{i=0}^n a_i$.
\end{proof}

\begin{cor} \label{bigGL} The $n \rightarrow \infty$ limit of the random
variables $\lambda_{\phi}$ with the uniform distribution on $GL(n,q)$ is
$P_{\frac{1}{q^{im_{\phi}}},
\frac{1}{q^{(i-1)m_{\phi}}},0,\frac{1}{q^{m_{\phi}}}}$.
\end{cor}

\begin{proof}
	Apply Lemma \ref{bign} to the final equation in the proof of
Theorem $\ref{main}$.
\end{proof}

	It is worth remarking that one possible motivation for a result
like Theorem \ref{main} comes from the theory of the symmetric groups, in
particular the ``Polya Cycle Index''. Let $a_i(\pi)$ be the number of
$i$-cycles of a permutation $\pi$. Using the fact that there are
$\frac{n!}{\prod_i a_i! i^{a_i}}$ elements of $S_n$ with $a_i$ $i$-cycles,
one proves that:

\[ \sum_{n=0}^{\infty} \frac{(1-u) u^n}{n!} \sum_{\pi \in S_n} \prod_i
x_i^{a_i(\pi)} = \prod_{m=1}^{\infty} e^{\frac{u^m}{m}(x_m-1)} \]

	This last equation has the following probabilistic
interpretation. Fix $u$ such that $0<u<1$ and pick the size of the
symmetric group with chance of size $n$ equal to $(1-u)u^n$. Next choose
$\pi$ uniformly in that $S_n$. Then the random variables $a_i(\pi)$ are
independent Poisson $\frac{u^i}{i}$. Combining the above equation with
Lemma \ref{bign} shows that for any $i< \infty$, the joint distribution of
$(a_1(\pi),\cdots,a_i(\pi))$ converges to independent (Poisson($1$),
$\cdots$, Poisson($\frac{1}{i}$)) as $n \rightarrow \infty$. The Poisson
distribution, naturally arising in the symmetric groups, is of fundamental
mathematical importance; it is reasonable to expect the distributions
$P_{\frac{u^{m_{\phi}}}{q^{im_{\phi}}},
\frac{1}{q^{(i-1)m_{\phi}}},0,\frac{1}{q^{m_{\phi}}}}$, naturally arising in the
general linear groups, to be of equal importance.

	Let us now consider briefly analogs of Theorem \ref{main} for the
other finite classical groups (proofs appear in Fulman \cite{fulthesis}).

\begin{enumerate}

\item {\bf Unitary Groups} The unitary group $U(n,q)$ can be defined as the
subgroup of $GL(n,q^2)$ preserving a non-degenerate skew-linear form, for
instance $<\vec{x},\vec{y}> = \sum_{i=1}^n x_i y_i^q$.
	
	Wall \cite{Wal} found that the conjugacy classes of the unitary
groups $U(n,q)$ have a nice combinatorial description analogous to rational
canonical form for the general linear groups. Given a polynomial $\phi$
with coefficients in $F_{q^2}$ and non-vanishing constant term, define a
polynomial $\tilde{\phi}$ by:

\[ \tilde{\phi} = \frac{z^{m_{\phi}} \phi^q(\frac{1}{z})}{[\phi(0)]^q} \]

	where $\phi^q$ raises each coefficient of $\phi$ to the $q$th
power. Writing this out, a polynomial $\phi(z)=z^{m_{\phi}} +
\alpha_{m_{\phi}-1} z^{m_{\phi}-1} + \cdots + \alpha_1 z + \alpha_0$ with
$\alpha_0 \neq 0$ is sent to $\tilde{\phi}(z)= z^{m_{\phi}} +
(\frac{\alpha_1}{\alpha_0})^q z^{m_{\phi}-1}+ \cdots +
(\frac{\alpha_{m_{\phi}-1}} {\alpha_0})^qz + (\frac{1}
{\alpha_0})^q$. Fulman \cite{fulthesis} proves that all irreducible
polynomials such that $\phi=\tilde{\phi}$ have odd degree.

	Wall $\cite{Wal}$ proves that the conjugacy classes of the
unitary group correspond to the following combinatorial data. As was
the case with $GL(n,q^2)$, an element $\alpha \in U(n,q)$ associates
to each monic, non-constant, irreducible polynomial $\phi$ over
$F_{q^2}$ a partition $\lambda_{\phi}$ of some non-negative integer
$|\lambda_{\phi}|$ by means of rational canonical form. The
restrictions necessary for the data $\lambda_{\phi}$ to represent a
conjugacy class are:

\begin{enumerate}
\item $|\lambda_z|=0$
\item $\lambda_{\phi}=\lambda_{\tilde{\phi}}$
\item $\sum_{\phi} |\lambda_{\phi}|m_{\phi}=n$
\end{enumerate}

	Random partitions $\lambda_{\phi}$ can be defined exactly as in the
case of the general linear groups. Fix $u$ such that $0<u<1$. Pick the size
with probability of size $n$ equal to $(1-u)u^n$. Next pick $\alpha$
uniformly in $U(n,q)$ and take the partition $\lambda_{\phi}(\alpha)$
corresponding to $\phi$ in Wall's description of the conjugacy class of
$\alpha$ in $U(n,q)$.

	Fulman \cite{fulthesis} uses Wall's conjugacy class size formula
and the fact that all polynomials invariant under $\tilde{}$ have odd
degree to prove the following analog of Theorem \ref{main}.

\begin{theorem} If $\phi = \tilde{\phi}$, then $\lambda_{\phi}$ has
distribution $P_{\frac{(-u)^{m_{\phi}}} {(-q)^{im_{\phi}}},\frac{1}
{(-q)^{(i-1)m_{\phi}}},0 ,\frac{1}{(-q)^{m_{\phi}}}}$. If $\phi \neq
\tilde{\phi}$, then $\lambda_{\phi}=\lambda_{\tilde{\phi}}$ have
distribution $P_{\frac{u^{2m_{\phi}}}
{q^{2im_{\phi}}},\frac{1}{q^{2(i-1)m_{\phi}}} ,0,\frac{1}
{q^{2m_{\phi}}}}$. These random partitions are independent and as with
$GL$, the case $u=1$ corresponds to the $n \rightarrow \infty$ limit.
\end{theorem}

\item {\bf Symplectic Groups} Assume for simplicity that the characteristic
of $F_q$ is not equal to 2. The symplectic group $Sp(2n,q)$ can be defined
as the subgroup of $GL(2n,q)$ preserving a non-degenerate alternating form
on $F_q$, for instance $<\vec{x},\vec{y}> = \sum_{i=1}^n (x_{2i-1}y_{2i} -
x_{2i}y_{2i-1})$.

	Given a polynomial $\phi$ with coefficients in $F_q$ and
non-vanishing constant term, define a polynomial $\bar{\phi}$ by:

\[ \bar{\phi} = \frac{z^{m_{\phi}} \phi^q(\frac{1}{z})}{[\phi(0)]^q} \]

	where $\phi^q$ raises each coefficient of $\phi$ to the $q$th
power. Explicitly, a polynomial $\phi(z)=z^{m_{\phi}} + \alpha_{m_{\phi}-1}
z^{m_{\phi}-1} + \cdots + \alpha_1 z + \alpha_0$ with $\alpha_0 \neq 0$ is
sent to $\bar{\phi}(z)= z^{m_{\phi}} + (\frac{\alpha_1}{\alpha_0})^q
z^{m_{\phi}-1}+ \cdots + (\frac{\alpha_{m_{\phi}-1}} {\alpha_0})^qz +
(\frac{1}{\alpha_0})^q$. (The notation $\bar{\phi}$ breaks from Wall
$\cite{Wal}$, in which $\tilde{\phi}$ was used, but these maps are
different. Namely $\tilde{}$ is defined on polynomials with coefficients in
$F_q$, but $\bar{}$ is defined on polynomials with coefficients in
$F_{q^2}$). Fulman \cite{fulthesis} showed that all irreducible polynomials
such that $\phi=\bar{\phi}$ have even degree, except for the polynomials $z
\pm 1$.

	Wall $\cite{Wal}$ proved that a conjugacy class of $Sp(2n,q)$
corresponds to the following data. To each monic, non-constant, irreducible
polynomial $\phi \neq z \pm 1$ associate a partition $\lambda_{\phi}$ of
some non-negative integer $|\lambda_{\phi}|$. To $\phi$ equal to $z-1$ or
$z+1$ associate a symplectic signed partition $\lambda_{\phi}^{\pm}$, by
which is meant a partition of some natural number $|\lambda_{\phi}^{\pm}|$
such that the odd parts have even multiplicity, together with a choice of
sign for the set of parts of size $i$ for each even $i>0$.

\begin{center}
Example of a Symplectic Signed Partition
\end{center}
	
\[ \begin{array}{c c c c c c c c}
	&  . & . & . & . & .  \\
	&  . & . & . & . & .  \\
	+  & . & . & . & . &  \\
	&  . & . & . &  &  \\
	&  . & . & . &  &   \\
	-  & . & . & & &   \\
	 & . & . & & &  
	  \end{array}  \]

	Here the $+$ corresponds to the parts of size 4 and the $-$
corresponds to the parts of size 2. This data represents a conjugacy
class of $Sp(2n,q)$ if and only if:

\begin{enumerate}

\item $|\lambda_{z}|=0$
\item $\lambda_{\phi}=\lambda_{\bar{\phi}}$
\item $\sum_{\phi=z \pm 1} |\lambda_{\phi}^{\pm}| + \sum_{\phi \neq z \pm 1} |\lambda_{\phi}| m_{\phi}=2n$

\end{enumerate}

	The symplectic groups can be used to define measures on partitions
$\lambda_{\phi}$ and symplectic signed partitions $\lambda^{\pm}_{z \pm 1}$
as follows. Fix $u$ so that $0<u<1$ and pick the dimension with probability
of dimension $2n$ equal to $(1-u^2)u^{2n}$. Then pick $\alpha$ uniformly in
$Sp(2n,q)$ and let $\lambda_{\phi}$ and $\lambda^{\pm}_{z \pm 1}$ be the
data corresponding to the conjugacy class of $\alpha$.

	Fulman \cite{fulthesis} uses Wall's conjugacy class size formula
and the fact that all polynomials other than $z \pm 1$ which are invariant
under $\bar{}$ have even degree to prove the following result.

\begin{theorem} If $\phi = \tilde{\phi}$ and $\phi \neq z \pm 1$, then
$\lambda_{\phi}$ has distribution $P_{\frac{(-u)^{m_{\phi}}}
{(-q)^{im_{\phi}}},\frac{1} {(-q)^{(i-1)m_{\phi}}},0
,\frac{1}{(-q)^{m_{\phi}}}}$. If $\phi \neq \tilde{\phi}$, then
$\lambda_{\phi}=\lambda_{\tilde{\phi}}$ have distribution
$P_{\frac{u^{2m_{\phi}}} {q^{2im_{\phi}}},\frac{1}{q^{2(i-1)m_{\phi}}}
,0,\frac{1} {q^{2m_{\phi}}}}$. These random partitions are independent and
as with $GL$, the case $u=1$ corresponds to the $n \rightarrow \infty$
limit.  \end{theorem}

	The distribution of the symplectic signed partitions
$\lambda^{\pm}_{z \pm 1}$ is more elusive (see Fulman \cite{fulthesis}) for
some results.
	
\item {\bf Orthogonal Groups} For simplicity assume that the characteristic
of $F_q$ is not equal to 2. The orthogonal groups can be defined as
subgroups of $GL(n,q)$ preserving a non-degenerate symmetric bilinear
form. For $n=2l+1$ odd, there are two such forms up to isomorphism, with
inner product matrices $A$ and $\delta A$, where $\delta$ is a non-square
in $F_q$ and $A$ is equal to:

\[ \left( \begin{array}{c c c}
		1 & 0 & 0 \\
		0 & 0_l & I_l \\
		0 & I_l & 0_l
	  \end{array} \right) \]

	Denote the corresponding orthogonal groups by $O^+(2l+1,q)$
and $O^-(2l+1,q)$. This distinction will be useful, even though these
groups are isomorphic.

	For $n=2l$ even, there are again two non-degenerate symmetric
bilinear forms up to isomorphism with inner product matrices:

\[  \left( \begin{array}{c c}
		0_l & I_l \\
		I_l & 0_l
	  \end{array} \right) \] 

\[    \left( \begin{array}{c c c c}
		0_{l-1} & I_{l-1} & 0 & 0 \\
		I_{l-1} & 0_{l-1} & 0 & 0 \\
		0 & 0 & 1 & 0 \\
		0 & 0 & 0 & -\delta
	  \end{array} \right)	\]

	where $\delta$ is a non-square in $F_q$. Denote the
corresponding orthogonal groups by $O^+(2l,q)$ and $O^-(2l,q)$.

	To describe the conjugacy classes of the finite orthogonal groups,
it is necessary to use the notion of the Witt type of a non-degenerate
quadratic form, as in Chapter 9 of Bourbaki $\cite{Bou}$. Call a
non-degenerate form $N$ null if the vector space $V$ on which it acts can
be written as a direct sum of 2 totally isotropic subspaces (a totally
isotropic space is one on which the inner product vanishes
identically). Define two non-degenerate quadratic forms $Q'$ and $Q$ to be
equivalent if $Q'$ is isomorphic to the direct sum of $Q$ and a null
$N$. The Witt type of $Q$ is the equivalence class of $Q$ under this
equivalence relation. There are 4 Witt types over $F_q$, which Wall denotes
by ${\bf 0},{\bf 1}, {\bf \delta}, {\bf \omega}$, corresponding to the
forms $0,x^2,\delta x^2,x^2 - \delta y^2$ where $\delta$ is a fixed
non-square of $F_q$. These 4 Witt types form a ring, but only the additive
structure is relevant here. The sum of two Witt types with representatives
$Q_1,Q_2$ on $V_1,V_2$ is the equivalence class of $Q_1+Q_2$ on
$V_1+V_2$. It is easy to see that the four orthogonal groups
$O^+(2n+1,q),O^-(2n+1,q), O^+(2n,q),O^-(2n,q)$ arise from forms $Q$ of Witt
types ${\bf 1},{\bf \delta},{\bf 0},{\bf \omega}$ respectively.

	Consider the following combinatorial data. To each monic,
non-constant, irreducible polynomial $\phi \neq z \pm 1$ associate a
partition $\lambda_{\phi}$ of some non-negative integer
$|\lambda_{\phi}|$. To $\phi$ equal to $z-1$ or $z+1$ associate an
orthogonal signed partition $\lambda_{\phi}^{\pm}$, by which is meant
a partition of some natural number $|\lambda_{\phi}^{\pm}|$ such that
all even parts have even multiplicity, and all odd $i>0$ have a choice
of sign. For $\phi= z-1$ or $\phi = z+1$ and odd $i>0$, we denote by
$\Theta_i (\lambda_{\phi}^{\pm})$ the Witt type of the orthogonal
group on a vector space of dimension $m_i(\lambda_{\phi}^{\pm})$ and
sign the choice of sign for $i$.

\begin{center}
Example of an Orthogonal Signed Partition
\end{center}
	
\[ \begin{array}{c c c c c}
	&  . & . & . & .   \\
	&  . & . & . & .   \\
	-  & . & . & . &   \\
	&  . & . &  &   \\
	&  . & . &  &    \\
	+  & . &  & &   \\
	 & . &  & &   
	  \end{array}  \]

	Here the $-$ corresponds to the part of size 3 and the $+$
corresponds to the parts of size 1. 

	The following theorem, though not stated there, is implicit in the
discussion on pages 38-40 of Wall $\cite{Wal}$. The polynomial $\bar{\phi}$
is defined as for the symplectic groups.

\begin{theorem} \label{almostWall} The data $\lambda^{\pm}_{z-1},
\lambda^{\pm}_{z+1}, \lambda_{\phi}$ represents a conjugacy class of
some orthogonal group if:

\begin{enumerate}
\item $|\lambda_{z}|=0$
\item $\lambda_{\phi}=\lambda_{\bar{\phi}}$
\item $\sum_{\phi=z \pm 1} |\lambda_{\phi}^{\pm}| + \sum_{\phi \neq z \pm
1} |\lambda_{\phi}| m_{\phi}=n$
\end{enumerate}
	
	In this case, the data represents the conjugacy class of exactly
one orthogonal group $O(n,q)$, with sign determined by the condition that
the group arises as the stabilizer of a form of Witt type:

\[ \sum_{\phi=z \pm 1} \sum_{i \ odd} \Theta_i(\lambda_{\phi}^{\pm}) +
\sum_{\phi \neq z \pm 1} \sum_{i \geq 1} i m_i(\lambda_{\phi}) {\bf \omega} \]

\end{theorem}

	The definition of measures on partitions for the tower $O(n,q)$
differs from that of the other groups. For $0<u<1$, pick an integer $n$
with the probability of $n=0$ equal to $\frac{1-u}{1+u}$ and probability of
$n \geq 1$ equal to $\frac{1-u}{1+u} 2u^n$. If $n \geq 1$, choose
$O^+(n,q)$ or $O^-(n,q)$ with probability $\frac{1}{2}$ and then choose
within that group uniformly. This defines random orthogonal signed
partitions $\lambda^{\pm}_{z-1}, \lambda^{\pm}_{z+1}$ and random partitions
$\lambda_{\phi}$ for $\phi \neq z \pm 1$. If $\phi \neq z \pm 1$, the
random partitions $\lambda_{\phi}$ have the same distribution as for the
symplectic groups. The orthogonal signed partitions are again elusive.
\end{enumerate}

\section{Hall-Littlewood Polynomials: The Young Tableau Algorithm}
\label{Tableau}

	This section, as the previous, studies the measures $P_{x,y,q,t}$
with $q=0, y_i=t^{i-1}$, and $x=ut^i$. We also set $t=\frac{1}{q}$ where
$q$, different from the $q$ above, is the size of a finite field. Section
\ref{Classical} showed that this is the case relevant to the finite
classical groups. As will emerge, the algorithm in this section is quite different
from the simplified algorithm of Section \ref{SimpleAlg}, which works by
adding horizontal strips.

	Recall that a standard Young tableau $T$ of size $n$ is a partition
of $n$ with each box containing one of $\{1,\cdots,n\}$ such that each of
$\{1,\cdots,n\}$ appears exactly once and the numbers increase in each row
and column of $T$. For instance,

\[ \begin{array}{c c c c c}
		1 & 3 & 5 & 6 &   \\
		2 & 4 & 7 &  &    \\
		8 & 9 &  &  &    
	  \end{array}  \]

	is a standard Young tableau. We call the algorithm in this section
the Young Tableau Algorithm because numbering the boxes in the order in
which they are created gives a standard Young tableau. It is assumed that
$0<u<1$ and $q>1$.

\begin{center}
The Young Tableau Algorithm
\end{center}

\begin{description}

\item [Step 0] Start with $N=1$ and $\lambda$ the empty
partition. Also start with a collection of coins indexed by the
natural numbers, such that coin $i$ has probability $\frac{u}{q^i}$ of
heads and probability $1-\frac{u}{q^i}$ of tails.

\item [Step 1] Flip coin $N$.

\item [Step 2a] If coin $N$ comes up tails, leave $\lambda$ unchanged,
set $N=N+1$ and go to Step 1.

\item [Step 2b] If coin $N$ comes up heads, choose an integer $S>0$
according to the following rule. Set $S=1$ with probability $\frac
{q^{N-\lambda_1'}-1} {q^N-1}$. Set $S=s>1$ with probability
$\frac{q^{N-\lambda_s'}-q^{N-\lambda_{s-1}'}}{q^N-1}$. Then increase
the size of column $s$ of $\lambda$ by 1 and go to Step 1.

\end{description}

	Note that as with the previous algorithms, this algorithm
halts by the Borel-Cantelli lemmas.

	Let us now look at the same example as in Section
$\ref{SimpleAlg}$, so as to see that the Young Tableau Algorithm is quite
different from the simplified algorithm for the Hall-Littlewood
polynomials.

	So suppose we are at Step 1 with $\lambda$ equal to the
following partition:
	
\[ \begin{array}{c c c c}
		. & . & . & .  \\
		. & . &  &      \\
		. &  &  &    \\
		 & & &  
	  \end{array}  \]

	Suppose also that $N=4$ and that coin 4 had already come up
heads once, at which time we added to column 1, giving $\lambda$. Now
we flip coin 4 again and get heads, going to Step 2b. We add to column
$1$ with probability $\frac{q-1}{q^4-1}$, to column $2$ with
probability $\frac{q^2-q}{q^4-1}$, to column $3$ with probability
$\frac{q^3-q^2} {q^4-1}$, to column $4$ with probability $0$, and to
column $5$ with probability $\frac{q^4-q^3}{q^4-1}$. We then return to
Step 1.

	Note that there is a non-0 probability of adding to column 1,
and that the dots added during the toss of a given coin need not form
a horizontal strip. This contrasts sharply with the algorithm in
Section $\ref{SimpleAlg}$.

	We use the notation that $(x)_N=(1-x)(1-\frac{x}{q}) \cdots
(1-\frac{x}{q^{N-1}})$. Recall from Section $\ref{Background}$ that
$m_i(\lambda)$ is the number of parts of $\lambda$ of size $i$, that
$n(\lambda) = \sum_i (i-1) \lambda_i$, that $a(s)$ is the number of squares
in $\lambda$ to the east of $s$, and that $l(s)$ is the number of squares
in $\lambda$ to the south of $s$. Lemma $\ref{coinNform}$ gives a formula
for the truncated measure $P^N_{\frac{u}{q^i},\frac{1}
{q^{i-1}},0,\frac{1}{q}}$ in terms of this notation.

\begin{lemma} \label{coinNform}
$P^N_{\frac{u}{q^i},\frac{1}{q^{i-1}},0,\frac{1}{q}} (\lambda) = 0$ if
$\lambda$ has more than $N$ parts. Otherwise:

\[ P^N_{\frac{u}{q^i},\frac{1}{q^{i-1}},0,\frac{1}{q}} (\lambda) =
\frac{u^{|\lambda|} (\frac{u}{q})_N
(\frac{1}{q})_N}{(\frac{1}{q})_{N-\lambda_1'}} \prod_{i \geq 1}
\frac{1}{q^{(\lambda_i')^2}(\frac{1}{q})_{m_i(\lambda)}} \]
 \end{lemma}

\begin{proof}
	The first statement is clear from the definition of the measure
$P^N_{\frac{u}{q^i},\frac{1}{q^{i-1}},0,\frac{1}{q}}$ in Section
$\ref{Defining}$. The second equality can be deduced from the definition of
$P^N_{\frac{u}{q^i},\frac{1}{q^{i-1}},0,\frac{1}{q}}$ and the Principal
Specialization Formula as follows:

\begin{eqnarray*}
P^N_{\frac{u}{q^i},\frac{1}{q^{i-1}},0,\frac{1}{q}} (\lambda) & = & \frac{P_{\lambda}(\frac{u}{q},\cdots,\frac{u}{q^N},0,\cdots;0,\frac{1}{q}) P_{\lambda}(\frac{1}{q^{i-1}};0,\frac{1}{q}) b_{\lambda}(0,t)}{\prod (\frac{u}{q},\cdots,\frac{u}{q^N},0,\cdots,\frac{1}{q^{i-1}};0,\frac{1}{q})}\\
& = & [\prod_{i=1}^N (1-\frac{u}{q^i})] P_{\lambda}(\frac{u}{q},\cdots,\frac{u}{q^N},0,\cdots;0,\frac{1}{q}) P_{\lambda}(\frac{1}{q^{i-1}};0,\frac{1}{q}) b_{\lambda}(0,t)\\
& = & \frac{\prod_{i=1}^N (1-\frac{u}{q^i}) (1-\frac{1}{q^i})}{\prod_{i=1}^{N-\lambda_1'} (1-\frac{1}{q^i})} \frac{P_{\lambda}(\frac{u}{q},\cdots,\frac{u}{q^N},0,\cdots;0,\frac{1}{q})}{q^{n(\lambda)}}\\
& = & \frac{u^{|\lambda|} (\frac{u}{q})_N (\frac{1}{q})_N}{(\frac{1}{q})_{N-\lambda_1'}}  \frac{P_{\lambda}(\frac{1}{q},\cdots,\frac{1}{q^N},0,\cdots;0,\frac{1}{q})}{q^{n(\lambda)}}\\
& = & \frac{u^{|\lambda|} (\frac{u}{q})_N (\frac{1}{q})_N}{(\frac{1}{q})_{N-\lambda_1'}} \frac{1}{q^{|\lambda|+2n(\lambda)}} \prod_{s \in \lambda: a(s)=0} \frac{1}{1-\frac{1}{q^{l(s)+1}}}\\
& = & \frac{u^{|\lambda|} (\frac{u}{q})_N
(\frac{1}{q})_N}{(\frac{1}{q})_{N-\lambda_1'}} \prod_{i \geq 1}
\frac{1}{q^{(\lambda_i')^2}(\frac{1}{q})_{m_i(\lambda)}}
\end{eqnarray*}

	where the last equality uses the fact that $n(\lambda) =
\sum_i {\lambda_i' \choose 2}$.
\end{proof}

\begin{theorem} \label{TableauAlg} The chance that the Young Tableau
algorithm yields $\lambda$ at the end of interval $N$ is
$P^N_{\frac{u}{q^i},\frac{1}{q^{i-1}},0,\frac{1}{q}}(\lambda)$.
\end{theorem}

\begin{proof}
	The theorem is clear if $N< \lambda_1'$ for then
$P^N_{\frac{u}{q^i},\frac{1}{q^{i-1}},0,\frac{1}{q}}(\lambda)=0$, and
Step 2b does not permit the number of parts of the partition to exceed
the number of the coin being tossed at any stage in the algorithm.

	For the case $N \geq \lambda_1'$, use induction on
$|\lambda|+N$. The base case is that $\lambda$ is the empty
partition. This means that coins $1,2,\cdots,N$ all came up tails on
their first tosses, which occurs with probability
$(\frac{u}{q})_N$. So the base case checks.

	Let $s_1 \leq s_2 \leq \cdots \leq s_k$ be the columns of
$\lambda$ with the property that changing $\lambda$ by decreasing the
size of one of these columns by 1 gives a partition
$\lambda^{s_i}$. It then suffices to check that the claimed formula
for $P^N_{\frac{u} {q^i},\frac{1}{q^{i-1}},0,\frac{1}{q}}(\lambda)$
satisfies the equation:

\begin{eqnarray*}
P^N_{\frac{u}{q^i},\frac{1}{q^{i-1}},0,\frac{1}{q}}(\lambda) & = & (1-\frac{u}{q^N})
P^{N-1}_{\frac{u}{q^i},\frac{1}{q^{i-1}},0,\frac{1}{q}}(\lambda) + \frac{u}{q^N}
\frac{q^{N-\lambda_1'}-1} {q^N-1} P^N_{\frac{u}{q^i},\frac{1}{q^{i-1}},0,\frac{1}{q}}(\lambda^{1})\\
& & + \sum_{s_i>1} \frac{u}{q^N} \frac{q^{N-\lambda_{s_i}'+1}- q^{N-\lambda_{s_i-1}'}}{q^N-1} P^N_{\frac{u}{q^i},\frac{1}{q^{i-1}},0,\frac{1}{q}}(\lambda^{s_i})
\end{eqnarray*}

	This equation is based on the following logic. Suppose that when
coin $N$ came up tails, the algorithm gave the partition $\lambda$. If coin
$N$ came up tails on its first toss, then we must have had $\lambda$ when
coin $N-1$ came up tails. Otherwise, for each $s_i$ we add the probability
that ``The algorithm gave the partition $\lambda^{s_i}$ on the penultimate
toss of coin $N$ and the partition $\lambda$ on the last toss of coin
$N$''. It is not hard to see that this probability is equal to the
probability of getting $\lambda^{s_i}$ on the final toss of coin $N$,
multiplied by the chance of a heads on coin $N$ which then gives the
partition $\lambda$ from $\lambda^{s_i}$.

	We divide both sides of this equation by $P^N_{\frac{u}{q^i}
,\frac{1}{q^{i-1}},0,\frac{1}{q}}(\lambda)$ and show that the terms on
the right-hand side sum to 1. First consider the terms with
$s_i>1$. Induction gives that:

\begin{eqnarray*}
& & \sum_{s_i>1} \frac{u}{q^N} \frac{q^{N-\lambda_{s_i}'+1}-
q^{N-\lambda_{s_{i-1}}'}}{q^N-1}
\frac{P^N_{\frac{u}{q^i},\frac{1}{q^{i-1}},0,\frac{1}{q}}(\lambda^{s_i})}{P^N_{\frac{u}{q^i},\frac{1}{q^{i-1}},0,\frac{1}{q}}(\lambda)}\\
& = & \sum_{s_i>1} \frac{q^{-\lambda_{s_i}'+1}- q^{-\lambda_{s_i-1}'}}{q^N-1}
\frac{q^{{\lambda_{s_i}' \choose 2}}
(\frac{1}{q})_{\lambda_{s_{i-1}}'-\lambda_{s_i}'}
(\frac{1}{q})_{\lambda_{s_i}' -
\lambda_{s_{i+1}}'}}{q^{{\lambda_{s_i}'-1 \choose 2}}
(\frac{1}{q})_{\lambda_{s_{i-1}}'-\lambda_{s_i}'+1}
(\frac{1}{q})_{\lambda_{s_i}' - \lambda_{s_{i+1}}'-1}}\\
& = & \sum_{s_i>1} \frac{q^{-\lambda_{s_i}'+1}-
q^{-\lambda_{s_{i-1}}'}}{q^N-1} q^{2 \lambda_{s_i}' -1}
\frac{(1-\frac{1}{q^{\lambda_{s_i}'-\lambda_{s_{i+1}}'}})}{(1-\frac{1}{q^{\lambda_{s_{i-1}}'-\lambda_{s_i}'+1}})}\\
& = & \sum_{s_i>1} \frac{q^{\lambda_{s_i}'}-q^{\lambda_{s_{i+1}}'}}{q^N-1}\\
& = &  \frac{q^{\lambda_2'}-1}{q^N-1}
\end{eqnarray*}

	Next consider the term coming from $P^{N-1}_{\frac{u}
{q^i},\frac{1}{q^{i-1}},0,\frac{1}{q}}(\lambda)$. If $N=\lambda_1'$, then
$\lambda_1' > N-1$, so by what we have proven
$P^{N-1}_{\frac{u}{q^i},\frac{1}{q^{i-1}} ,0,\frac{1}{q}} (\lambda)$ is
0. Otherwise,

\begin{eqnarray*}
(1-\frac{u}{q^N}) \frac{P^{N-1}_{\frac{u}{q^i},\frac{1}{q^{i-1}},0,\frac{1}{q}}(\lambda)}{
P^N_{\frac{u}{q^i},\frac{1}{q^{i-1}},0,\frac{1}{q}}(\lambda)} & = & (1-\frac{u}{q^N})
\frac{(\frac{u}{q})_{N-1} (\frac{1}{q})_{N-1}
(\frac{1}{q})_{N-\lambda_1'}}{(\frac{u}{q})_N (\frac{1}{q})_N
(\frac{1}{q})_{N-\lambda_1'-1}}\\
& = & \frac{(1-\frac{1}{q^{N-\lambda_1'}})}{(1-\frac{1}{q^N})}\\
& = & \frac{q^N-q^{\lambda_1'}}{q^N-1}
\end{eqnarray*}

	So this term always contributes $\frac{q^N-q^{\lambda_1'}}
{q^N-1}$.

	Finally, consider the term coming from
$P^N_{\frac{u}{q^i},\frac{1}{q^{i-1}},0,\frac{1}{q}}
(\lambda^{1})$. This vanishes if $\lambda_1=\lambda_2$ since then
$\lambda^{1}$ is not a partition. Otherwise,

\begin{eqnarray*}
\frac{u}{q^N} \frac{q^{N-\lambda_{s_1}'+1}-1}{q^N-1}
\frac{P^N_{\frac{u}{q^i},\frac{1}{q^{i-1}},0,\frac{1}{q}}(\lambda^{s_1})}{P^N_{\frac{u}{q^i},\frac{1}{q^{i-1}},0,\frac{1}{q}}(\lambda)}
& = & \frac{q^{-\lambda_{s_1}'+1}- q^{-N}}{q^N-1} \frac{(\frac{1}{q})_{N-\lambda_1'}}{(\frac{1}{q})_{N-\lambda_1'+1}}
\frac{q^{{\lambda_{s_1}' \choose 2}}
(\frac{1}{q})_{\lambda_{1}' -
\lambda_{2}'}}{q^{{\lambda_{1}'-1 \choose 2}}
(\frac{1}{q})_{\lambda_{1}' - \lambda_{2}'-1}}\\
& = & \frac{q^{-\lambda_1'+1}- q^{-N}}{q^N-1}
\frac{1}{(1-\frac{1}{q^{N-\lambda_1'+1}})} q^{2 \lambda_1' -1}
(1-\frac{1}{q^{\lambda_1'-\lambda_2'}})\\
& = & \frac{q^{\lambda_1'}-q^{\lambda_2'}}{q^N-1}
\end{eqnarray*}

	So in all cases this term contributes $\frac{q^{\lambda_1'} -
q^{\lambda_2'}}{q^N-1}$.

	Adding up the three terms completes the proof.
\end{proof}

	As an example of Lemma $\ref{coinNform}$ and Theorem
$\ref{TableauAlg}$, suppose that $N=4$ and $\lambda$ is the partition:
	
\[ \begin{array}{c c}
		. & . \\
		. & \\
		. & 
	  \end{array}  \]

	Then the chance that the Young Tableau Algorithm gives the
partition $\lambda$ when coin $4$ comes up tails is:

\[ \frac{u^4(1-\frac{u}{q})(1-\frac{u}{q^2})(1-\frac{u}{q^3})
(1-\frac{u}{q^4})(1-\frac{1}{q^3})(1-\frac{1}{q^4})}{q^{10}(1-\frac{1}{q})^2}
\]

\section{Hall-Littlewood Polynomials: Weights on the Young Lattice} \label{Weights}

	In this section $T$ denotes a standard Young tableau and $\lambda$
denotes the partition corresponding to $T$. Let $|T|$ be the size of
$T$. As explained in Section $\ref{Tableau}$, the Young Tableau algorithm
constructs a standard Young tableau, and thus defines a measure on the set
of all standard Young tableaux.

	Let $P_{\frac{u}{q^i},\frac{1}{q^{i-1}},0,\frac{1}{q}}(T)$ be
the chance that the Young Tableau algorithm of Section $\ref{Tableau}$
outputs $T$, and let $P^N_{\frac{u}{q^i} ,\frac{1}{q^{i-1}}
,0,\frac{1}{q}}(T)$ be the chance that it outputs $T$ when coin $N$
comes up tails.

	We also introduce the following notation. Let $T_{(i,j)}$ be
the entry in the $(i,j)$ position of $T$ (recall that $i$ is the row
number and $j$ the column number). For $j \geq 2$, let $A_{(i,j)}$ be
the number of entries $(i',j-1)$ such that $T_{(i',j-1)} <
T_{(i,j)}$. Let $B_{(i,j)}$ be the number of entries $(i',1)$ such
that $T_{(i',1)} < T_{(i,j)}$. For instance the tableau:

\[ \begin{array}{c c c c c}
		1 & 3 & 5 & 6 &   \\
		2 & 4 & 7 &  &    \\
		8 & 9 &  &  &    
	  \end{array}  \]

	has $T_{(1,3)}=5$. Also $A_{(1,3)}=2$ because there are 2
entries in column $3-1=2$ which are less than 5 (namely 3 and
4). Finally, $B_{(1,3)}=2$ because there are 2 entries in column $1$
which are less than 5 (namely 1 and 2).

	There is a simple formula for $P^N_{\frac{u}{q^i},
\frac{1}{q^{i-1}},0,\frac{1}{q}}(T)$ in terms of this data.

\begin{theorem} \label{Lattice}
$P^N_{\frac{u}{q^i},\frac{1}{q^{i-1}},0,\frac{1}{q}}(T)=0$ if $T$ has
greater than $N$ parts. Otherwise:

\[ P^N_{\frac{u}{q^i},\frac{1}{q^{i-1}},0,\frac{1}{q}}(T) =
\frac{u^{|T|}}{|GL(\lambda_1',q)|} \frac{\prod_{r=1}^N
(1-\frac{u}{q^r}) (1-\frac{1}{q^r})}{\prod_{r=1}^{N-\lambda_1'}
(1-\frac{1}{q^r})} \prod_{(i,j) \in \lambda \atop j \geq 2}
\frac{q^{1-i}-q^{-A_{(i,j)}}}{q^{B_{(i,j)}}-1} \]

\end{theorem}

\begin{proof}
	The case where $T$ has more than $N$ parts is proven as in Theorem
$\ref{TableauAlg}$.
	
	The case $\lambda_1' \leq N$ is proven by induction on
$|T|+N$. If $|T|+N=1$, then $T$ is the empty tableau and $N=1$. This
means that coin 1 in the Tableau algorithm came up tails on the first
toss, which happens with probability $1-\frac{u}{q}$. So the base case
checks.

	For the induction step, there are two cases. The first case
is that the largest entry in $T$ occurs in column $s>1$. Removing the
largest entry from $T$ gives a tableaux $T^{s}$. We have the
equation:

\[ P^N_{\frac{u}{q^i},\frac{1}{q^{i-1}},0,\frac{1}{q}}(T) =
(1-\frac{u}{q^N}) P^{N-1}_{\frac{u}{q^i},
\frac{1}{q^{i-1}},0,\frac{1}{q}}(T) + \frac{u}{q^N}
\frac{q^{N-\lambda_s'+1}-q^{N-\lambda_{s-1}'}}{q^N-1}
P_{\frac{u}{q^i},\frac{1}{q^{i-1}},0,\frac{1}{q}}^N(T^{s}) \]

	The two terms in this equation correspond to the whether or
not $T$ was completed at time $N$. We divide both sides of the
equation by $P^N_{\frac{u}{q^i},\frac{1}{q^{i-1}},0,\frac{1}{q}}(T)$,
substitute in the conjectured formula, and show that it satisfies this
recurrence. The two terms on the right hand side then give:

\[ \frac{q^N-q^{\lambda_1'}}{q^N-1} + \frac{q^
{-\lambda_{s}'+1}-q^{-\lambda_{s-1}'}} {q^N-1} \frac{1}
{\frac{q^{-\lambda_{s}'+1}-q^{-\lambda_{s-1}'}}{q^{\lambda_1'}-1}} = 1
\]

	The other case is that the largest entry of $T$ occurs in
column 1. We then have the equation:

\[ P_{\frac{u}{q^i},\frac{1}{q^{i-1}},0,\frac{1}{q}}^N(T) =
(1-\frac{u}{q^N})P_{\frac{u}{q^i},\frac{1}{q^{i-1}},0,\frac{1}{q}}^{N-1}(T)
+ \frac{u}{q^N} \frac{q^{N-\lambda_1'}-1}{q^N-1}
P_{\frac{u}{q^i},\frac{1}{q^{i-1}},0,\frac{1}{q}}^N(T^{1}) \]

	As in the previous case, we divide both sides of the equation
by $P_{\frac{u}{q^i},\frac{1}{q^{i-1}},0,\frac{1}{q}}^N(T)$,
substitute in the conjectured formula, and show that it satisfies this
recurrence. The two terms on the right hand side then give:

\[ \frac{q^N-q^{\lambda_1'}}{q^N-1} + \frac{1}{q^N}
\frac{q^{N-\lambda_1'+1}-1}{q^N-1} \frac{1}
{(1-\frac{1}{q^{N-\lambda_1'+1}})} \frac{|GL(\lambda_1',q)|}
{|GL(\lambda_1'-1,q)|} = 1 \]

	This completes the induction, and the proof of the theorem.
\end{proof}

	For instance, Theorem $\ref{Lattice}$ says that if

\[ S = \frac{u^4(1-\frac{u}{q})(1-\frac{u}{q^2})(1-\frac{u}{q^3})(1-\frac{u}{q^4}) (1-\frac{1}{q^2})(1-\frac{1}{q^3})(1-\frac{1}{q^4})}{|GL(3,q)|} \]

	then the chances that the Young Tableau Algorithm gives the
following tableaux:

\[ \begin{array}{c c}
		1 & 2 \\
		3 & \\
		4 & 
	  \end{array}  \]

\[ \begin{array}{c c}
		1 & 3 \\
		2 & \\
		4 & 
	  \end{array}  \]

\[ \begin{array}{c c}
		1 & 4 \\
		2 & \\
		3 & 
	  \end{array}  \]

	 when coin $4$ comes up tails are $\frac{S}{q}$,
$\frac{S}{q^2}$, and $\frac{S}{q^3}$ respectively. Note that the sum
of these probabilities is:

\[ \frac{u^4(1-\frac{u}{q})(1-\frac{u}{q^2})(1-\frac{u}{q^3})
(1-\frac{u}{q^4})(1-\frac{1}{q^3})(1-\frac{1}{q^4})}{q^{10}(1-\frac{1}{q})^2}
\]
	
	As must be the case and as was proved at the end of Section
$\ref{Tableau}$, this quantity is also equal to the chance that the
Young Tableau Algorithm gives the partition:
	
\[ \begin{array}{c c}
		. & . \\
		. & \\
		. & 
	  \end{array}  \]

	An interesting object in combinatorics is the Young lattice. The
elements of this lattice are all partitions of all numbers. An edge is
drawn between partitions $\lambda$ and $\Lambda$ if $\Lambda$ is obtained
from $\lambda$ by adding one box. Note that a standard Young tableau $T$ of
shape $\lambda$ is equivalent to a path in the Young lattice from the empty
partition to $\lambda$. This equivalence is given by growing the partition
$\lambda$ by adding boxes in the order $1,\cdots,n$ in the positions
determined by $T$. For instance the tableau:

\[ \begin{array}{c c c}
		1 & 3 & 4 \\
		2 &  &    \\    
	  \end{array}  \]

	corresponds to the path:

\[ \begin{array}{c c c c c c c c c c c c c c c c c c c c}
		 & & & & . & & & & . &&&& .&. & & && .&.&.                          \\
		 \emptyset & & \rightarrow &&  && \rightarrow && . && \rightarrow && . &&& \rightarrow && .& & \\    
	  \end{array}  \]

	The measure $P_{\frac{u}{q^i},\frac{1}{q^{i-1}} ,0,\frac{1}{q}}(T)$
on standard Young tableaux has the following description in terms of
weights on the Young lattice.

\begin{cor} \label{Weight} Put weights $m_{\lambda,\Lambda}$ on
the Young lattice according to the rules:

\begin{enumerate}

\item $m_{\lambda,\Lambda} = \frac{u}{q^{\lambda_1'}(q^{\lambda_1'+1}-1)}$ if
$\Lambda$ is obtained from $\lambda$ by adding a box to column 1

\item $m_{\lambda,\Lambda} = \frac{u(q^{-\lambda_s'}-q^{-
\lambda_{s-1}'})}{q^{\lambda_1'}-1}$ if $\Lambda$ is obtained from
$\lambda$ by adding a box to column $s>1$

\end{enumerate}

	Then the chance that the Tableau algorithm produces $T$ is
equal to:

\[ \prod_{r=1}^{\infty} (1-\frac{u}{q^r}) \prod_{i=0}^{|T|-1}
m_{\gamma_i,\gamma_{i+1}} \]

	where the $\gamma_i$ are the partitions in the path along the
Young lattice which corresponds to the tableau $T$.
\end{cor}

\begin{proof}
	This follows by letting $N \rightarrow \infty$ in Theorem
$\ref{Lattice}$ and the fact that $T$ corresponds to a unique path in
the Young lattice.
\end{proof}

	The following remarks may be of interest.

\begin{enumerate}

\item Note that the total weight out of the empty partition is
$\frac{u}{q-1}$ and that the total weight out of any other partition
$\lambda$ is:

\begin{eqnarray*}
\frac{u}{q^{\lambda_1'}(q^{\lambda_1'+1}-1)} + \sum_{i \geq 2}
\frac{u(q^{-\lambda_s'}-q^{- \lambda_{s-1}'})}{q^{\lambda_1'}-1} & = & \frac{u}{q^{\lambda_1'}(q^{\lambda_1'+1}-1)} + \frac{u}{q^{\lambda_1'}}\\
& = & \frac{uq}{q^{\lambda_1'+1}-1}\\
& < & 1
\end{eqnarray*}

	Since the sum of the weights out of a partition $\lambda$ to a
larger partition $\Lambda$ is less than 1, the weights can also be
viewed as transition probabilities, provided that one allows for
halting.

\item Note that the Young Tableau Algorithm of Section $\ref{Tableau}$ for
growing $\lambda_{\phi}$ according to the group theoretic measures of
Section \ref{Tableau} does not carry over to unitary case if
$\phi=\tilde{\phi}$, since then some of the probabilities involved would be
negative. The description in terms of weights on the Young lattice in
Corollary $\ref{Weight}$, however, does extend to the unitary groups. The
weight formula should be altered as follows. In the case
$\phi=\tilde{\phi}$ one replaces the variables $(u,q)$ by $(-u,-q)$, and in
the case $\phi \neq \tilde{\phi}$ one replaces the variables $(u,q)$ by
$(u^2,q^2)$.

\item Some applications of the results of this and the preceding section
toward proving group theoretic theorems can be found in the companion paper
by Fulman \cite{fulalgorithm}.

\end{enumerate}

\section{Schur Functions: A $q$-analog of the Plancherel Measure of the
Symmetric Group} \label{Plancherel}

	To begin, let us recall the definition of the Plancherel measure of
the symmetric group. This is a measure on the partitions $\lambda$ of size
$n$. Letting $h(s)=a(s)+l(s)+1$ be the hook-length of $s \in \lambda$, the
Plancherel measure assigns to $\lambda$ the probability $\frac{n!}{\prod_{s
\in \lambda} h(s)^2}$. Kerov and Vershik \cite{Ker}, \cite{Ver1},
\cite{Ver2} have studied Plancherel measure extensively. The connection with the
representation theory of the symmetric group is that the irreducible
representations of $S_n$ can be parameterized by partitions $\lambda$ of $n$ such
that the representation corresponding to $\lambda$ has dimension
$\frac{n!}{\prod_{s \in \lambda} h(s)}$ (see pages 53-96 of Sagan $\cite{Sagan}$).

	Plancherel measure has another description. Robinson and Schensted
found a bijection from the symmetric group to the set of pairs $(P,Q)$ of
standard Young tableau of the same shape (see pages 97-101 of
Sagan \cite{Sagan} for details). Call the shape associated to $\pi$ under the
Robinson-Schensted correspondence $\lambda(\pi)$. Then $\lambda(\pi)$ has
Plancherel measure if $\pi$ is chosen uniformly from the symmetric group. This
follows from the fact that the dimension of the
irreducible representation of $S_n$ corresponding to the partition
$\lambda$ is the number of standard tableaux of shape $\lambda$.

	Let us now see how the measures $P_{x,y,q,t}$ of Section
\ref{Defining} lead to a $q$-analog of Plancherel measure. This section
studies the specialization $x_i=t^i, y_i=t^{i-1}, q=t$. We then set
$t=\frac{1}{q}$, where this $q$ is the size of a finite field.
	Lemma $\ref{formula}$ gives a formula for the measure
$P^N_{\frac{u}{q^i},\frac{1}{q^{i-1}},\frac{1}{q},\frac{1}{q}}$. We use the
notation that $(x)_N=(1-x)(1-\frac{x}{q}) \cdots (1-\frac{x}{q^{N-1}})$. Let
$c(s)=a'(s)-l'(s)$ denote the content of $s \in \lambda$ (here
$l'(s),l(s),a(s)$, and $a'(s)$ are the number of squares in $\lambda$ to
the north, south, east, and west of $s$ respectively).

\begin{lemma} \label{formula}

\[ P^N_{\frac{1}{q^i},\frac{1}{q^{i-1}} ,
\frac{1}{q},\frac{1}{q}}(\lambda) = [\prod_{r=1}^N
\prod_{t=0}^{\infty} (1-\frac{1}{q^{r+t}})]
\frac{1}{q^{2n(\lambda)+|\lambda|}} \prod_{s \in \lambda}
\frac{1-\frac{1}{q^{N+c(s)}}}{(1-\frac{1}{q^{h(s)}})^2} \]

\end{lemma}

\begin{proof}
	This can be deduced from the definition of the measure
$P^N_{\frac{1}{q^i},\frac{1}{q^{i-1}},\frac{1}{q},\frac{1}{q}}$ and
the Principal Specialization Formula as follows:

\begin{eqnarray*}
P^N_{\frac{1}{q^i},\frac{1}{q^{i-1}},\frac{1}{q},\frac{1}{q}}(\lambda) & = &
\frac{P_{\lambda}(\frac{1}{q},\cdots,\frac{1}{q^N},0,\cdots;\frac{1}{q},\frac{1}{q})
P_{\lambda}(\frac{1}{q^{i-1}};\frac{1}{q},\frac{1}{q})
b_{\lambda}(\frac{1}{q},\frac{1}{q})}{\prod
(\frac{1}{q},\cdots,\frac{1}{q^i},0,\cdots,\frac{1}{q^{i-1}})}\\
& = & [\prod_{r=1}^N \prod_{t=0}^{\infty} (1-\frac{1}{q^{r+t}})] \frac{1}{q^{n(\lambda)}} \prod_{s \in \lambda} \frac{1}{(1-\frac{1}{q^{h(s)}})} P_{\lambda}(\frac{1}{q},\cdots,\frac{1}{q^N},0,\cdots;\frac{1}{q},\frac{1}{q})\\
& = &  [\prod_{r=1}^N \prod_{t=0}^{\infty} (1-\frac{1}{q^{r+t}})]
\frac{1}{q^{2n(\lambda)+|\lambda|}} \prod_{s \in \lambda}
\frac{1-\frac{1}{q^{N+c(s)}}}{(1-\frac{1}{q^{h(s)}})^2}
\end{eqnarray*}

\end{proof}

	Renormalizing the measure $P_{\frac{1}{q^i},\frac{1}{q^{i-1}}
, \frac{1}{q},\frac{1}{q}}$ to live on partitions of size $n$ will
give a $q$-analog of the Plancherel measure. To this end, we introduce
polynomials $J_n(q)$. First define $J_{\lambda}(q)$ by:

\[ J_{\lambda}(q) = \frac{q^{|\lambda|^2-|\lambda|-2n(\lambda)}
[(\frac{1}{q})_{|\lambda|}]^2} {\prod_{s \in \lambda}
(1-\frac{1}{q^{h(s)}})^2} \]

	The measure $P_{\frac{1}{q^i}, \frac{1}{q^{i-1}} ,
\frac{1}{q},\frac{1}{q}}$ can then be written as:

\[ P_{\frac{1}{q^i},\frac{1}{q^{i-1}} ,
\frac{1}{q},\frac{1}{q}}(\lambda) = [\prod_{r=1}^{\infty}
\prod_{t=0}^{\infty} (1-\frac{1}{q^{r+t}})]
\frac{J_{\lambda}(q)}{q^{|\lambda|^2} (1-\frac{1}{q})^2 \cdots
(1-\frac{1}{q^{|\lambda|}})^2} \]

	It is not clear that the $J_{\lambda}(q)$ are
polynomials in $q$, but this will turn out to be true. Define $J_n(q)
= \sum_{\lambda: |\lambda|=n} J_{\lambda}(q)$ and $J_0(q)=1$.	Proposition
$\ref{condit}$, which follows immediately from the definitions in this section,
explains why one might be interested in the polynomials $J_{\lambda}(q)$ and
$J_n(q)$.

\begin{prop} \label{condit} Under the measure $P_{\frac{1}{q^i},\frac{1}{q^{i-1}},\frac{1}{q},\frac{1}{q}}$, the conditional probability of $\lambda$ given that $|\lambda|=n$ is equal to $\frac{J_{\lambda}(q)}{J_n(q)}$.
\end{prop}

	Lemma \ref{hook} is an easy exercise from page 11 of Macdonald \cite{Mac} and
will be useful.
	
\begin{lemma} \label{hook}

\[\sum_{s \in \lambda} h(s) =  n(\lambda)+n(\lambda')+|\lambda|\]

\end{lemma}

	It is possible to relate the polynomials $J_{\lambda}(q)$ to
the Kostka-Foulkes polynomials $K_{\lambda}(q)$ (sometimes
denoted $K_{\lambda (1^n)}(q)$). The Kostka-Foulkes polynomials are defined as:

\[ K_{\lambda}(q) = \frac{q^{n(\lambda)} [|\lambda|]!}{\prod_{s \in \lambda}
[h(s)]} \]

	where $[n]=1+q+\cdots+q^{i-1}$, the $q$-analog of the number
$n$. One can also check from Chapter 4 of Macdonald $\cite{Mac}$
that $K_{\lambda'}(q)$ is the degree of the unipotent representation of
$GL(n,q)$ corresponding to the partition $\lambda'$.

	Proposition $\ref{Relate}$ connects the $J_{\lambda}(q)$ to
the Kostka-Foulkes polynomials.

\begin{prop} \label{Relate} $J_{\lambda}(q) = [K_{\lambda'}(q)]^2$
\end{prop}

\begin{proof}
	Using Lemma $\ref{hook}$, observe that:

\begin{eqnarray*}
J_{\lambda}(q) & = & \frac{q^{|\lambda|^2-|\lambda|-2n(\lambda)}
[(\frac{1}{q})_{|\lambda|}]^2} {\prod_{s \in \lambda}
(1-\frac{1}{q^{h(s)}})^2}\\
& = & q^{|\lambda|^2-|\lambda|-2n(\lambda)} \frac{q^{2 \sum_{s \in \lambda} h(s)}} {\prod_{s \in \lambda} (q^{h(s)}-1)^2} \frac{\prod_{i=1}^{|\lambda|} (q^i-1)^2}{q^{|\lambda|^2+|\lambda|}}\\
& = & q^{2 \sum_{s \in \lambda} h(s) - 2|\lambda| -2n(\lambda)} (\frac{[|\lambda|]!}{\prod_{s \in \lambda} [h(s)]})^2\\
& = & q^{2n(\lambda')} (\frac{[|\lambda'|]!}{\prod_{s \in \lambda'} [h(s)]})^2\\
& = & K_{\lambda'}(q)^2
\end{eqnarray*}

\end{proof}

	Theorem $\ref{property}$ gives some properties of the $J_n(q)$. By
the remark before Proposition $\ref{Relate}$, $J_n(q)$ is the sum of the
squares of the degrees of the irreducible unipotent representations of
$GL(n,q)$. Recall that $[u^n] f(u)$ means the coefficient of $u^n$ in
$f(u)$.

\begin{theorem} \label{property}

\begin{enumerate}

\item $J_n(q)$ is a symmetric polynomial of degree $2 {n \choose 2}$
which has non-negative integer coefficients and satisfies $J_n(1)=n!$.

\item $\frac{J_n(q)}{q^{n^2}(1-\frac{1}{q})^2 \cdots
(1-\frac{1}{q^n})^2} = [u^n] \frac{1}{\prod_{r=1}^{\infty}
\prod_{s=0}^{\infty} (1-\frac{u}{q^{r+s}})}$

\end{enumerate}

\end{theorem}

\begin{proof}
	Proposition $\ref{Relate}$ shows that $J_n(q)$ is a polynomial
with non-negative integer coefficients. Note by Lemma $\ref{hook}$
that:

\begin{eqnarray*}
deg (J_{\lambda}) & = & 2deg(K_{\lambda'})\\
& = & 2[n(\lambda ')+ {|\lambda|+1 \choose 2} - \sum_{s \in \lambda'} h(s)]\\
& = & 2 {|\lambda| \choose 2} - 2n(\lambda)
\end{eqnarray*}

	Thus $J_{\lambda}$ has degree $2 {|\lambda| \choose 2}$ for
$\lambda=(|\lambda|)$ and smaller degree for all other $\lambda$. So
$J_n(q)$ has degree $2 {n \choose 2}$. Symmetry means that $J_n(q) = q^{2{n
\choose 2}} J_n(\frac{1}{q})$. In fact $J_{\lambda}(q)+J_{\lambda'}(q)$ satisfies
this property, by Lemma $\ref{hook}$.

	To see that $J_n(1)=n!$, observe that:

\begin{eqnarray*}
J_n(1) & = & \sum_{\lambda \vdash n} [K_{\lambda'}(1)]^2\\
& = & \sum_{\lambda \vdash n} [K_{\lambda}(1)]^2\\
& = & \sum_{\lambda \vdash n} [\frac{n!}{\prod_{s \in \lambda} h(x)}]^2\\
& = & n!
\end{eqnarray*}

	For the second part of the theorem, it is useful to consider the
measure $P_{\frac{u}{q^i},\frac{1}{q^{i-1}} ,\frac{1}{q}
,\frac{1}{q}}$. Arguing as in Lemma $\ref{formula}$ shows that:

\[P_{\frac{u}{q^i},\frac{1}{q^{i-1}},\frac{1}{q},\frac{1}{q}}(\lambda) =
[\prod_{r=1}^{\infty} \prod_{t=0}^{\infty} (1-\frac{u}{q^{r+t}})]
\frac{u^{|\lambda|} J_{\lambda}(q)}{q^{|\lambda|^2} (1-\frac{1}{q})^2
\cdots (1-\frac{1}{q^{|\lambda|}})^2} \]

	The fact that this is a measure means that:

\[ \sum_{n=1}^{\infty} \frac{u^n J_n(q)}{q^{n^2}(1-\frac{1}{q^2})
\cdots (1-\frac{1}{q^n})^2} = \frac{1}{ [\prod_{r=0}^{\infty}
\prod_{s=1}^{\infty} (1-\frac{u}{q^{r+s}})]} \]

	Taking coefficients of $u^n$ on both sides proves the second
part.
\end{proof}
  
	Corollary $\ref{con}$ of Theorem $\ref{property}$ shows that
conditioning the measure $P_{\frac{1}{q^i}, \frac{1}{q^{i-1}} ,
\frac{1}{q},\frac{1}{q}}$ on $|\lambda|=n$ gives a $q$-analog of the
Plancherel measure on partitions of size $n$.

\begin{cor} $\label{con}$ The conditional probability of $\lambda$ given
that $|\lambda|=n$ under the measure $P_{\frac{1}{q^i},\frac{1}{q^{i-1}} ,
\frac{1}{q},\frac{1}{q}}$ reduces to the Plancherel measure of the
symmetric group when one sets $q=1$.
\end{cor}

\begin{proof}
	Proposition $\ref{condit}$ shows that the conditional
probability is $\frac{J_{\lambda}(q)}{J_n(q)}$. The result follows
from the definition of $J_{\lambda}(q)$, and the fact that
$J_n(1)=n!$, which is part of the first statement of Theorem
$\ref{property}$.
\end{proof}

	The following observations show that this $q$-analog of Plancherel
measure has properties similar to the Plancherel measure of the symmetric
group.

\begin{enumerate}

\item By Proposition \ref{Relate} and the remark before it, our $q$-analog
of Plancherel measure assigns a probability to $\lambda$ which is
proportional to the square of the degree of the unipotent representation of
$GL(n,q)$ parameterized by $\lambda'$, the transpose partition. This is in
direct analogy with the Plancherel measure of the symmetric group, which
assigns a probability to $\lambda$ which is proportional to the square of
the degree of the irreducible representation of $S_n$ parameterized by
$\lambda'$.

\item The description of the Plancherel measure of the symmetric group in
terms of the Robinson-Schensted correspondence carries over to the above
$q$-analog of Plancherel measure. To state this precisely recall that the
major index of a permutation $\pi \in S_n$ is defined by:

\[ maj(\pi) = \sum_{i: 1 \leq i \leq n-1 \atop \ \  \pi(i)>\pi(i+1)} i \]

\begin{theorem} Choose $\pi \in S_n$ with probability proportional to
$q^{maj(\pi)+maj(\pi^{-1})}$. Then $\lambda(\pi)'$, the transpose of the
partition associated to $\pi$ through the Robinson-Schensted
correspondence, has the $q$-analog of Plancherel measure defined in
Corollary \ref{con}.
\end{theorem}

\begin{proof}
	Define the major index of a standard Young tableau as the sum of
the entries $i$ such that $i+1$ is in a row below that of $i$. Reasoning
similar to that of page 243 of Macdonald \cite{Mac} shows that

\[ K_{\lambda'}(q) = \sum_{T \in SYT(\lambda')} q^{maj(T)} \]

	where the sum is over all standard Young tableaux of shape
$\lambda'$. 

	From the way the Robinson-Schensted correspondence works (pages
97-101 of Sagan \cite{Sagan}), one sees that if $\pi$ corresponds to the
pair $(P,Q)$, then $maj(\pi)=maj(Q)$. It is also known (Theorem 3.86 of
Sagan \cite{Sagan}) that if $\pi$ corresponds to the pair $(P,Q)$, then
$\pi^{-1}$ corresponds to the pair $(Q,P)$.

	Proposition \ref{Relate} and the Robinson-Schensted correspondence
thus give that:

\begin{eqnarray*}
J_{\lambda}(q) & = & [\sum_{T \in SYT(\lambda')} q^{maj(T)}]^2\\
& = & \sum_{(P,Q) \in \{SYT(\lambda') \times SYT(\lambda')\}} q^{maj(P)}
q^{maj(Q)}\\
& = & \sum_{\pi \in S_n: \lambda(\pi)=\lambda'} q^{maj(\pi)+maj(\pi^{-1})}
\end{eqnarray*}

\end{proof}

\end{enumerate}

\section{Schur Functions: A Comparison with Kerov's $q$-analogs of
Plancherel Measure and the Hook Walk} \label{Kerov}

	Kerov $\cite{Kerq}$ has a $q$-analog of Plancherel measure which
comes from the Schur functions. His $q$-analog of Plancherel measure is defined
implicitly by means of a probabilistic algorithm called the $q$ hook walk. This
walk starts with the empty partition, and adds a box at a time. The
partition $\lambda$ grows to $\Lambda$ (here $|\Lambda|=|\lambda|+1$)
with probability:

\[ \frac{q^{n(\Lambda)} \prod_{s \in \lambda} [h(s)]}{q^{n(\lambda)}
\prod_{s \in \Lambda} [h(s)]} \]

	It can now be seen that Kerov's $q$-analog of Plancherel
measure is different from the $q$-analog introduced in Section
$\ref{Plancherel}$, because the partition

\[ \begin{array}{c c}
		. & .   
	  \end{array}  \]

	has mass $\frac{1}{q+1}$ under Kerov's $q$-analog of
Plancherel measure and mass $\frac{q^2}{q^2+1}$ under our $q$-analog
of Plancherel measure.

	Proposition $\ref{equiv}$ relates Kerov's $q$ hook walk to the
algorithm of Section $\ref{MainAlg}$.

\begin{prop} \label {equiv} Suppose that $n_N$ is equal to 1 for all $N$ in
Step 1 of the algorithm of Section $\ref{MainAlg}$ for picking from
$P_{\frac{1}{q^i},\frac{1}{q^{i-1}} ,\frac{1}{q},\frac{1}{q}}$. The growth
process on partitions this defines is exactly Kerov's $q$ hook walk.
\end{prop}

\begin{proof}
	Step 2 in the algorithm of Section $\ref{MainAlg}$
changes $\lambda$ to $\Lambda$ with probability:

\[ \frac{\phi_{\Lambda /
\lambda}(\frac{1}{q},\frac{1}{q})}{g_1(\frac{1}{q^{i-1}})}
\frac{P_{\Lambda}(\frac{1}{q^{i-1}};
\frac{1}{q},\frac{1}{q})}{P_{\lambda}(\frac{1}{q^{i-1}};\frac{1}{q},\frac{1}{q})}
\]

	The definition of $\phi_{\Lambda / \lambda}$ shows that
$\phi_{\Lambda / \lambda}(\frac{1}{q},\frac{1}{q})=1$. Corollary
$\ref{Addn}$ shows that $g_1=\frac{1}{1-\frac{1}{q}}$. The Principal
Specialization Formula shows that $P_{\lambda}(\frac{1}{q^{i-1}};
\frac{1}{q},\frac{1}{q})$ is equal to $\frac{1}{q^{n(\lambda)}}
\prod_{s \in \lambda} \frac{1}{1-q^{h(s)}}$. Combining these facts
proves that:

\[ \frac{\phi_{\Lambda /
\lambda}(\frac{1}{q},\frac{1}{q})}{g_1(\frac{1}{q^{i-1}})}
\frac{P_{\Lambda}(\frac{1}{q^{i-1}};
\frac{1}{q},\frac{1}{q})}{P_{\lambda}(\frac{1}{q^{i-1}};\frac{1}{q},\frac{1}{q})} =  \frac{q^{n(\Lambda)} \prod_{s \in \lambda} [h(s)]}{q^{n(\lambda)}
\prod_{s \in \Lambda} [h(s)]} \]

	as desired.
\end{proof}

\section{Suggestions for Future Research} \label{Suggestions}

	This section suggests some possibilities for future research.

\begin{enumerate} 

\item Read group theoretic information off of the probabilistic algorithms
of Sections \ref{Tableau} and \ref{Weights}. As was shown in Section
\ref{Classical}, these algorithms are related to the finite classical
groups. Some group theoretic applications of these algorithms are given in
the companion paper \cite{fulalgorithm}.

\item Develop probabilistic algorithms for picking from the measures
$\lambda_{z \pm 1}^{\pm}$ for the symplectic and orthogonal groups. These
will be more complicated than the algorithms for the general linear and
unitary groups, since there are size restrictions on the partitions (for
instance in the symplectic groups $|\lambda_{z \pm 1}^{\pm}|$ is always
even). Presumably one adds $1*2$ or $2*1$ tiles according to some
rules.

\item Persi Diaconis suggested the problem of implementing this paper's
algorithms in a computer program. The Young
Tableau Algorithm, for instance, involves flipping infinitely many coins. How can
this practical obstacle be overcome?

\item Study the shapes of partitions under the measures
$P_{x,y,q,t}(\lambda)$ for various specializations of the variables
$x,y,q,t$. For instance find generating functions for various functionals
of the partitions such as the number of parts, largest part, number of 1's,
etc. (A generating function for the size was found as Corollary
$\ref{size}$ of Section $\ref{MainAlg}$). It should also be possible to
extend work of Vershik $\cite{Ver1}, \cite{Ver2}$ which shows that random
partitions under measures such as the Plancherel measure have an asymptotic
limit shape.

\end{enumerate}

\section{Acknowledgments} The work here is taken from the author's Ph.D.
thesis, done under the guidance of Persi Diaconis at Harvard University. The
author thanks him for numerous ideas and suggestions. The author also thanks
Arkady Berenstein and A.N. Kirillov for their help. This research was done under
the generous 3-year support of the National Defense Science and Engineering
Graduate Fellowship (grant no. DAAH04-93-G-0270) and the support of the
Alfred P. Sloan Foundation Dissertation Fellowship.


\begin{thebibliography}{AAA}

\bibitem{And} Andrews, G., The theory of partitions. Encyclopedia of
Mathematics and its Applications, Vol. 2. Addison-Wesley Publishing Co.,
Reading, Mass.-London-Amsterdam, 1976. 

\bibitem{ArT} Arratia, R. and Tavare, S., The cycle structure of random
permutations. Ann. Probab. 20 (1992), no. 3, 1567-1591.

\bibitem{Bou} Bourbaki, N., Formes sesqulineaires et formes
quadratiques (Elements de Mathematique I, livre II), Hermann (Paris),
1959.

\bibitem{Ca1} Carter, R., Simple groups of lie type. John Wiley and
Sons, 1972.

\bibitem{Nu3} Celler, F., Leedham-Green, C., Murray, S., Niemeyer, A.,
and O'Brien, E.A., Generating random elements of a finite
group. Communications in algebra, 23 (13), (1995), 4931-4948.

\bibitem{Di2} Diaconis, P. and Shahshahani, M., On the eigenvalues of
random matrices, J. Appl. Prob. 31 (1994), 49-61.

\bibitem{FiH} Fine, N.J. and Herstein, I. N., The probability that a
matrix is nilpotent, Illinois J. Math. 2 (1958), 499-504.

\bibitem{Fri} Fristedt, B., The structure of random partitions of large
integers. Trans. Amer. Math. Soc. 337 (1993), no. 2, 703-735.

\bibitem{fulthesis} Fulman, J., Probability in the classical groups over
finite fields: symmetric functions, stochastic algorithms and cycle
indices, PhD Thesis, Harvard University, 1997.

\bibitem{fulalgorithm} Fulman, J., A probabilistic approach toward
the finite general linear and unitary groups, preprint.

\bibitem{Ger} Gerstenhaber, M., On the number of nilpotent matrices
with coefficients in a finite field, Illinois J. Math. 5 (1961), 330-333.

\bibitem{Go2} Goh, W. and Schmutz, E., A central limit theorem on
$GL_n(F_q)$, Preprint. Department of Math. Drexel University. 

\bibitem{Greene} Greene, C., Nijenhuis, A. and Wilf, H., A
probabilistic proof of a formula for the number of Young tableaux of a
given shape. Adv. in Math 31 (1979), no. 1, 104-109.

\bibitem{Greene2} Greene, C., Nijenhuis, A. and Wilf, H., Another
probabilistic method in the theory of Young
tableaux. J. Combin. Theory Series A, 37 (1984), 127-135.

\bibitem{Han} Hansen, J. and Schmutz, E., How random is the characteristic
polynomial of a random matrix? Math. Proc. Cambridge Philos. Soc. 114
(1993), no. 3, 507-515.

\bibitem{Her} Herstein, I.N., Topics in algebra. Second edition. Xerox
College Publishing, Lexington, Mass.-Toronto, Ont., 1975.

\bibitem{Kerov} Kerov, S.V., The boundary of Young lattice and random
Young tableaux. Formal power series and algebraic
combinatorics. DIMACS Ser. Discrete Math. Theoret. Comput. Sci. 24,
pg. 133-158.

\bibitem{Kerq} Kerov, S.V., A $q$-analog of the hook walk algorithm
for random Young tableaux. Journal of Algebraic Combinatorics 2
(1993), 383-396.

\bibitem{Ker} Kerov, S.V. and Vershik, A.M., Asymptotic behavior of
the Plancherel measure of the symmetric group and the limit form of
Young tableaux. Dokl. Akad. Nauk SSSR 233 (1977), no. 6, 1024-1027.

\bibitem{Kun} Kung, J., The cycle structure of a linear transformation over
a finite field, Linear Algebra Appl. 36 (1981), 141-155.

\bibitem{Mac} Macdonald, I.G., Symmetric functions and Hall
polynomials, Second Edition. Claredon Press, Oxford. 1995.

\bibitem{Meh} Mehta, M.L., Random matrices. Academic Press, San
Diego. (1991). 

\bibitem{Sagan} Sagan, B., The symmetric group: representations,
combinatorial algorithms, and symmetric functions. Wadsworth and
Brooks/Cole 1991.

\bibitem{St1} Stong, R., Some asymptotic results on finite vector
spaces, Advances in Applied Mathematics 9, 167-199 (1988).

\bibitem{St2} Stong, R., The average order of a matrix. Journal
of Combinatorial Theory, Series A. Vol. 64, No. 2, November 1993.

\bibitem{Ver1} Vershik, A.M., Asymptotic combinatorics and algebraic
analysis. Proceedings of the International Congress of Mathematicians,
Zurich 1994, 1384-1394.

\bibitem{Ver2} Vershik, A.M., Statistical mechanics of combinatorial
partitions, and their limit shapes. Functional Analysis and its
Applications, Vol. 30, No. 2, 1996, pg. 90-105.

\bibitem{Wal} Wall, G.E., On conjugacy classes in the unitary,
symplectic, and orthogonal groups, Journal of the Australian
Mathematical Society 3 (1963), 1-63.

\end{thebibliography}
\end{document}